 \documentclass[final,3p,times,twocolumn]{elsarticle}

\usepackage{amssymb}
\usepackage{amsmath}
\usepackage{xcolor}
\usepackage{algpseudocode}
\usepackage{algorithm2e}
\usepackage{float}

\usepackage{subcaption}
\usepackage{bm}
\usepackage{booktabs}
\usepackage{multirow}

\journal{Ocean Engineering}

\begin{document}
\setlength{\parskip}{0pt}
\begin{frontmatter}



\title{A multi-objective optimization framework \\
for reducing the impact of ship noise on marine mammals}

\author[1]{Akash Venkateshwaran}
\author[1]{Indu Kant Deo}
\author[2]{Jasmin Jelovica}
\author[1]{Rajeev K. Jaiman}


\affiliation[1]{organization={Department of Mechanical Engineering},
          addressline={The University of British Columbia},
           city={Vancouver},
           postcode={V6T 1Z4},
           state={BC},
           country={Canada}}
\affiliation[2]{organization={Department of Mechanical and Civil Engineering},
           addressline={The University of British Columbia},
           city={Vancouver},
           postcode={V6T 1Z4},
           state={BC},
           country={Canada}}



\begin{abstract}
The underwater radiated noise (URN) emanating from ships presents a significant threat to marine mammals, given their heavy reliance on hearing. The intensity of URN from ships is correlated to their speed, making speed reduction a crucial operational mitigation strategy. This paper presents a new multi-objective optimization framework to optimize the ship speed for effective URN mitigation without compromising fuel consumption. This framework addresses a fixed-path voyage scheduling problem, incorporating two objective functions namely, noise intensity levels and fuel consumption. The optimization is performed using the state-of-the-art non-dominated sorting genetic algorithm under voyage constraints. A 2D ocean acoustic environment, comprising randomly scattered marine mammals of diverse audiogram groups and realistic conditions, including sound speed profiles and bathymetry, is simulated. To estimate the objective functions, we consider empirical relations for fuel consumption and near-field noise modeling together with a ray-tracing approach for far-field noise propagation. The optimization problem is solved to determine the Pareto solutions and the trade-off solution. The effectiveness of the framework is demonstrated via practical case studies involving a large container ship. A comparative analysis illustrates the adaptability of the framework across different oceanic environments, affirming its potential as a robust tool for reducing the URN from shipping.


\end{abstract}




\end{frontmatter}


\section{Introduction}

\begin{figure*}[ht]
\centering
\includegraphics[width=15cm]{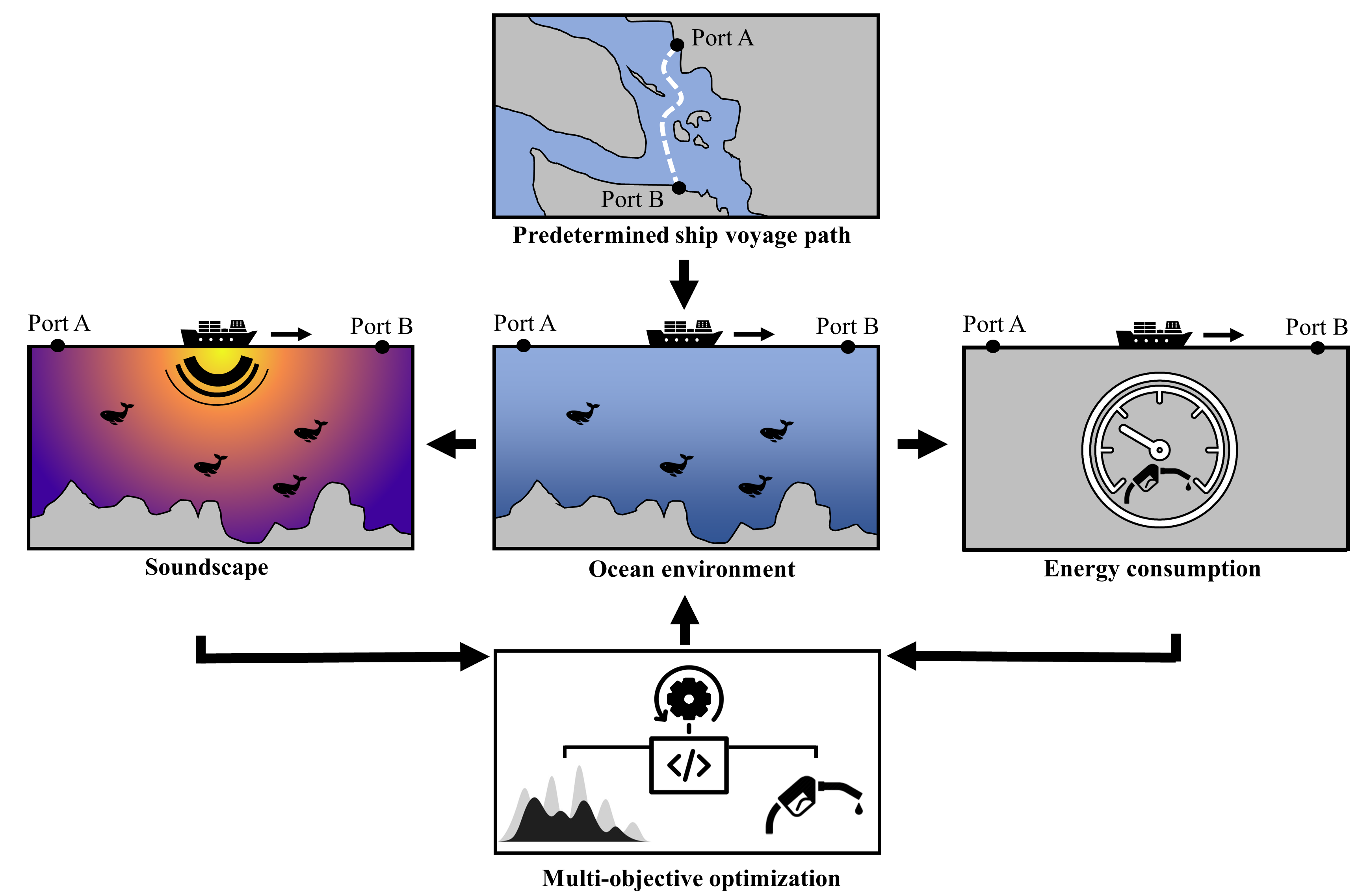}
\caption{A high-level overview of the proposed optimization framework.}
\label{fig_general_flow_chart}
\end{figure*}

Maritime transport has always been a major part of international
commercial transportation (more than 80\%), as it constitutes a highly cost-effective way of transferring large volumes of cargo between continents. As a result, there has been a dramatic upsurge in the volume of global seaborne cargo and is anticipated to increase over the coming years \cite{Kaplan2016}. The increase in maritime activity has prompted detrimental consequences, including increased emissions of greenhouse gases and chemical pollutants, incidents of marine mammal collisions, and a rise in underwater noise pollution \cite{Jagerbrand2019}. Particularly, the latter remains a significant concern to marine engineers, biologists and public policymakers.

Modern ships have become the primary contributors to anthropogenic noise in the oceans responsible for increasing the ambient noise levels at low frequencies (10-1k Hz) at a rate of 3 dB/decade up until the 1980s \cite{Peng2015, MiksisOlds2016,McDonald2006}. Consequently, marine mammals pose a particular concern, given their heavy reliance on hearing for essential activities such as foraging, orientation, predator detection, and communication at these frequencies. Underwater radiated noise (URN) emitted by vessels can directly interfere with their behavior, leading to disturbances in these activities \cite{Erbe2019, Richardson1995}. Numerous studies are underway, reporting the physiological and behavioral responses to anthropogenic noise \cite{Gomez2016, Wright2007, Southall2019}.

This increasing human acoustic footprint has led to collective mitigation measures including regulatory actions, spatiotemporal management strategies for noise sources, and advancements in vessel quieting technologies \cite{Chou2021}. Among these measures, voluntary speed reduction or vessel slowdowns have emerged as an effective means of mitigation and management \cite{Leaper2019, Joy2019}. This strategy reduces the acoustic footprint by 95\% following a 30\% speed reduction and concurrently decreases the overall ship strike rate by 50\% \cite{Leaper2019}. The Enhancing Cetacean Habitat and Observation Program (ECHO) projects, led by the Vancouver Fraser Port Authority, include seasonal vessel slowdown trials in Haro Strait, Boundary Pass, and Swiftsure Bank to reduce and evaluate URN in the feeding area of southern resident killer whales. The 2022 trial in Haro Strait and Boundary Pass had 93\% pilot-reported participation. The results reported a median reduction of broadband noise by 2.7 dB (equivalent to a 46\% reduction in sound intensity) within Haro Strait and a corresponding reduction of 2.8 dB (equivalent to a 48\% reduction in sound intensity) within Boundary Pass. Moreover, foraging conditions were improved with whales being detected 26\% of the days the slowdown trial was active \cite{ECHO_boundarypass}. Similarly, there was a 3.1 dB (equivalent to a 51\% reduction in sound intensity) reduction in broadband sound intensity associated with the 2022 trial in Swiftsure Bank \cite{ECHO_swiftsure}. These studies indicate that the intensity of the noise is positively correlated with the speed of the vessel, as a result, speed reduction is an effective mitigation measure.

This paper presents an adaptive operational strategy to mitigate URN from shipping, specifically on selecting the optimal vessel speed in a voyage scheduling problem such that the overall impact of noise on mammals is minimized. A comprehensive overview of ship speed optimization models and the fundamental principles of the ship voyage and routing problem is presented in \cite{Psaraftis2014}. Previous papers on ship voyage optimization include minimizing main engine fuel consumption, and the ship operating costs \cite{Psaraftis2013}; considering varying sea states and weather conditions \cite{Li2020, Tzortzis2021}; emphasizing greenhouse gas and emission control-based voyages \cite{Yu2021, Ma2020}; combinations of emission and cost reduction as multi-objective optimization problems (MOOPs) \cite{Ma2021, Wang2019}; and multiple ship routing problems \cite{Wen2017}. Given these papers, the main idea behind them is the direct relationship between the objectives (e.g., propulsion power, operating cost, or emission rate) and the ship's speed and optimizing it while satisfying the voyage constraints. Most of this research regarding ship voyage optimization lacks consideration of the impact of URN.

In this paper, we address this gap in the literature by developing a ship speed optimization framework that considers the URN from shipping combined with fuel consumption. A multi-objective optimization framework (MOOF) is developed where the URN emission and the fuel consumption are quantified as objective functions while satisfying constraints related to the voyage. Fig. \ref{fig_general_flow_chart} shows the framework broadly where a ship travels from Port A to Port B on a predetermined route. A multi-objective optimization algorithm optimizes the ship speed based on inputs from the soundscape and energy consumption. A former work presented a speed optimization problem in all-electric ships where URN levels were maintained below a threshold as a constraint \cite{Khatami2023}. The URN was represented as a function of the vessel's speed at the source level. Nevertheless, it's essential to acknowledge that the received noise levels by marine mammals can vary significantly from the source levels \cite{Stojanovic2009}. The received noise levels are dependent on several factors such as the propagation of sound, its frequency characteristics, and various environmental properties (such as bathymetry, sound speed profiles, and seabed properties). Using the information of the received levels and the hearing sensitivity of the mammals, one can precisely ascertain the possibility of species being affected at varying distances from the noise source \cite{Southall2019}. Hence, we herein model the propagation of URN from the ship to get the received noise level for the objective function. Fig. \ref{fig_general_flow_chart} illustrates a general workflow of the proposed optimization framework, including a ship as the source of noise, the distribution of acoustic intensity levels in the ocean environment, the inputs of the acoustic and energy consumption to the multi-objective optimization algorithm.

The acoustic propagation model estimates the acoustic field at the receiver location given the input environment condition. These models leverage our understanding of sound propagation in the ocean to calculate the acoustic field. The acoustic wave equation is employed to derive these physics-based models \cite{desanto1979theoretical}.
Closed-form analytical solutions to the wave equation are intractable for varying ocean environments \cite{oliveira2021underwater}. There are various approximate solutions to the wave equation, the majority of which can be classified into the following four groups: ray methods, normal modes, parabolic equations, and wave number integration \cite{jensen2011computational}. The ray method is particularly attractive for high-frequency ($>$ 500 Hz) \cite{etter2018underwater}, range-dependent problems whereby the normal mode
or the parabolic-equation models are not practical alternatives. However, the standard ray-tracing method produces certain artifacts. The Gaussian beam method associates with each ray a beam with a Gaussian intensity profile normal to the ray. In comparison to standard ray tracing, the Gaussian beam method has the advantage of being free of ray-tracing artifacts \cite{porter1987gaussian}. In this work, we employ the Gaussian beam tracing method implemented in the bellhop solver \cite{porter2011bellhop} to generate the acoustic transmission loss at the receiver locations.

The paper is structured as follows. In Section \ref{sec_prob_descript}, the description of the problem and underlying assumptions are given. Subsequently, Section \ref{sec_math_forml} presents the mathematical formulation of the optimization problem including the acoustic model, fuel consumption model, and formulation of the objective functions as well as the constraints. The case studies are demonstrated in Section \ref{sec_case_studies} followed by the results in Section \ref{sec_results_d}. Finally, the paper is concluded with final remarks and prospects for future work in Section \ref{sec_conclusion}.

\section{Problem description} \label{sec_prob_descript}

In this paper, we study a multi-objective voyage optimization problem of a ship traveling in a fixed route, from Port A to Port B. The two objective functions are the total intensity of the URN signals received by the mammals and the total fuel consumption of the voyage. Ideally, reducing the vessel speed will decrease URN levels, however, this will cause an increase in fuel rate consumption and the sailing time. Consequently, these two conflicting objective functions necessitate the application of multi-objective optimization where we specifically use a posterior method \cite{Silva2016}.

\begin{figure*}[ht]
\centering
\includegraphics[width=15cm]{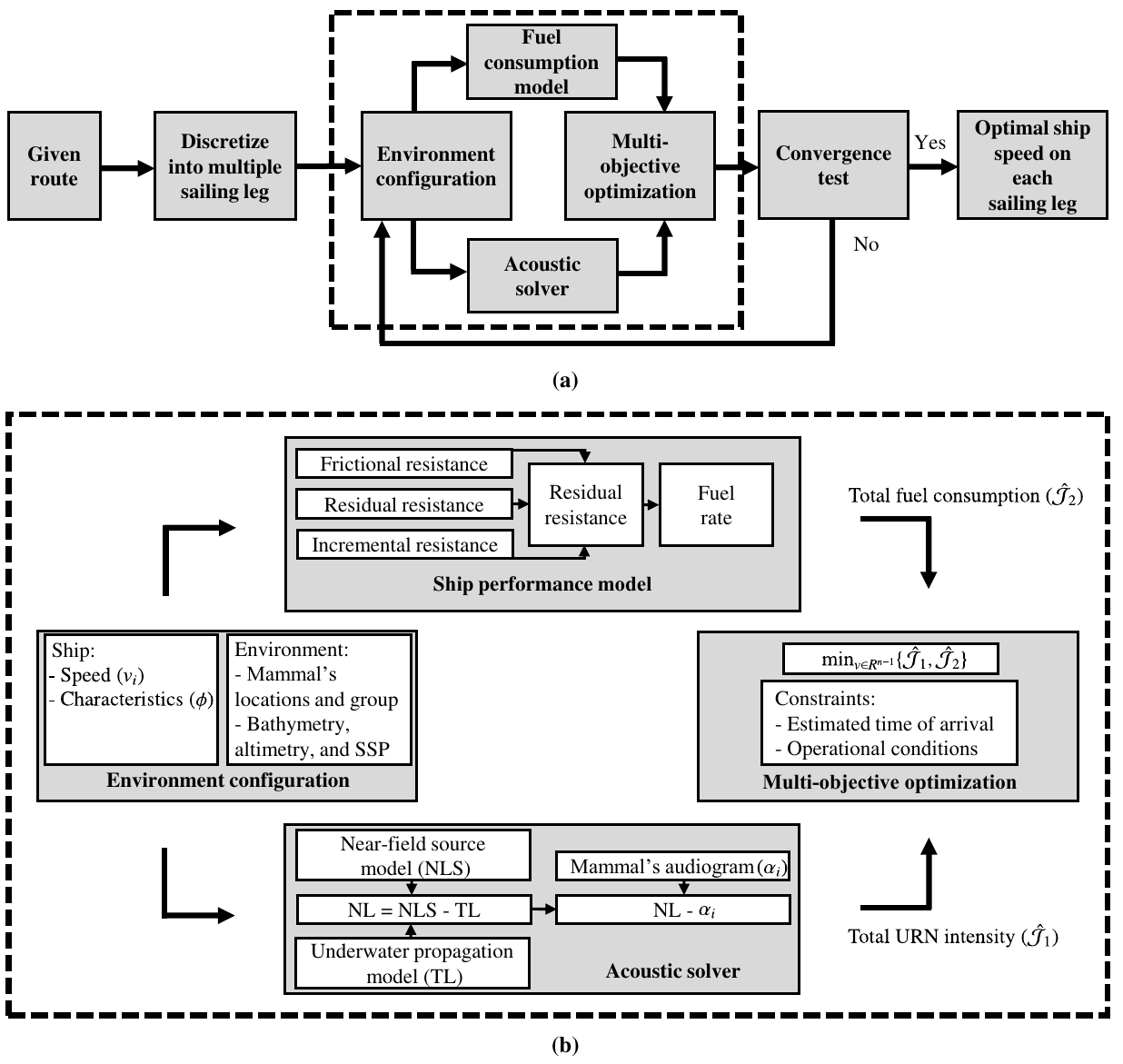}
\caption{Multi-objective optimization framework: (a) Main pipeline; (b) Modelling pipeline}
\label{fig_flow_chart}
\end{figure*}

Figure  \ref{fig_flow_chart}(a) illustrates the multi-objective optimization framework. The given route is discretized into multiple equidistant sections known as sailing legs, and the algorithm determines the optimal ship speed for each leg. Based on the route, a 2D ocean environment is formulated, encompassing the ship, randomly scattered mammals, and other environmental factors. By considering the operation conditions of the ship on a particular sailing leg, an empirical fuel consumption model calculates the total consumption for each sailing leg. Simultaneously, an acoustic model, incorporating information about the environment, simulates the ship's URN throughout the voyage. The calculated fuel consumption and the modeled noise levels are then fed as objective functions into an optimization algorithm which then determines the optimal vessel speed for each leg in compliance with the specified voyage constraints. 


Figure \ref{fig_flow_chart}(b) shows a comprehensive overview of calculations behind each of these components of the framework. The environmental configuration includes marine mammals distributed in fixed locations, as well as bathymetry, altimetry, and sound speed profile. Subsequently, we employ a near-field noise source model based on the environmental configuration to determine the noise source level at each sailing leg. In parallel, the acoustic wave equation is solved over a range of frequencies using an underwater propagation model to obtain the noise level received by each mammal. Additionally, from the environmental configuration, we use a ship performance model that calculates total resistance, which translates into fuel consumption rate. These two quantities are then defined as objective functions and solved using an optimization algorithm. Detailed mathematical descriptions are presented in the following sections.

This framework relies on several key assumptions, which are stated as follows:
\begin{enumerate}
\item The ship's speed remains constant throughout each sailing leg.
\item The payload or displacement of the ship remains constant throughout the entire voyage.
\item The environmental conditions, such as the sound speed profile, bathymetry, and altimetry, are predetermined and remain constant throughout the voyage.
\item The locations of marine mammals along the route are randomly scattered and fixed, and the associated audiogram group to which they belong is known.
\item The propagation of URN from the ship is simulated only when it is at a waypoint. Hence, the noise level received by the mammals is assumed to be constant along a given sailing leg.
\end{enumerate}

\begin{figure*}[ht]
\centering
\includegraphics[width=16cm]{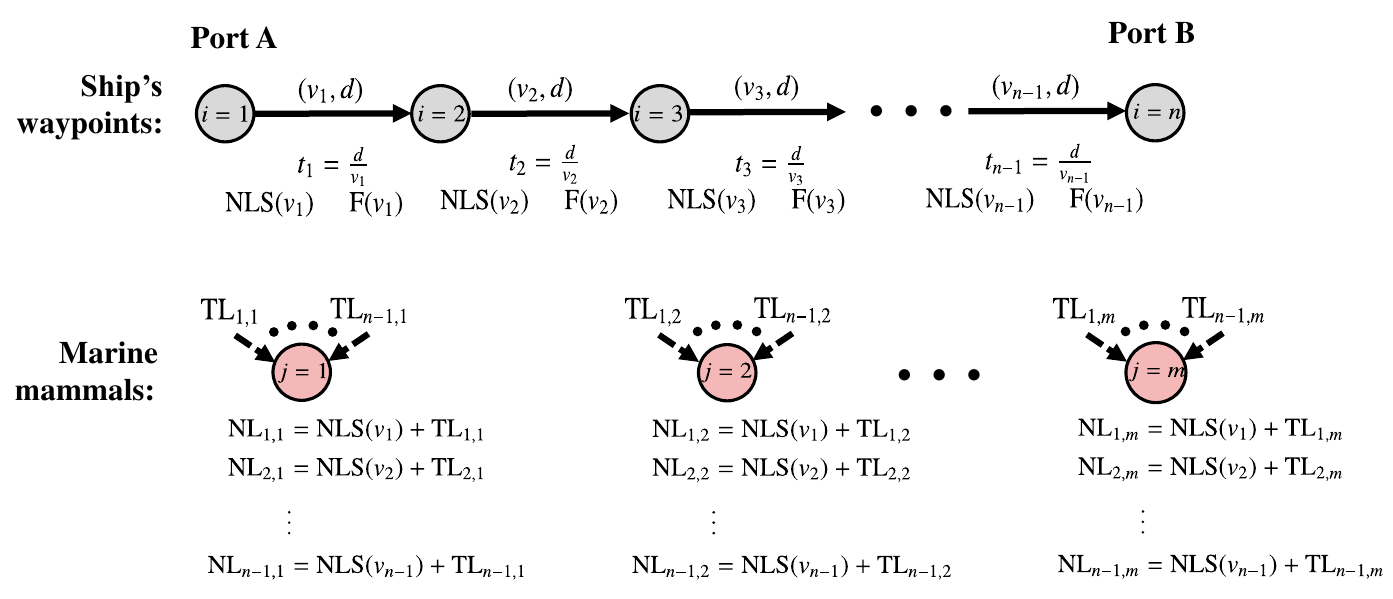}
    \caption{Ship-marine mammal network for estimating received noise levels by mammals at every position of the ship.}
\label{fig_net}
\end{figure*}

\section{Mathematical modelling} \label{sec_math_forml}

In this section, we introduce the underlying models, objective functions, and constraints pertinent to the optimization problem. The ship's route is segmented into multiple equidistant sailing legs as mentioned before. These legs are joined by two nodes or waypoints indexed by i = $\{1, 2, ..., n\}$. Each sailing leg has a specific velocity ($v_i$), i.e., the decision variable, and equal distance ($d$). Given the velocity and distance of a sailing leg, the time required to complete the leg ($t$), as well as the corresponding noise levels at the source ($\mathrm{NLS}(\cdot)$) and fuel consumption rate ($\mathrm{F}(\cdot)$) can be determined.

Let's consider a scenario where we have multiple marine mammals along the path. Each mammal is indexed by j = ${1, 2, ..., m}$, with a unique audiogram ($\alpha_j$) associated with it. Further details of the audiogram functions are provided in the upcoming section. It is important to note that as the ship progresses through each sailing leg, its velocity and position relative to the mammal constantly change, consequently the noise level received by the mammal. The received noise level ($\mathrm{NL}(\cdot)$) is determined by the difference between the noise level at the source ($\mathrm{NLS}(\cdot)$) and the transmission loss ($\mathrm{TL}(\cdot)$) of the noise signal. This is repeated for every sailing leg ($i$) and all the mammals ($j$) to compute the overall intensity of the noise signal, accounting for the entire voyage and expressed as an objective function. Similarly, the fuel consumption rate for each sailing leg is aggregated to determine the total fuel used for the entire voyage, which serves as another objective function.

Figure \ref{fig_net} illustrates the segmentation of the route and calculations within each sailing leg. To summarise, when the ship sails from waypoint $i$ to $i+1$ at velocity $v_i$, the fuel consumption rate, $F_i(\cdot)$, and the noise level received at each mammal location, $\mathrm{NL}_{i,j}(v_i) = \mathrm{NLS}(v_i) - \mathrm{TL}_{i,j}$, for all $j \in \{1,2,..., m \}$, are calculated. A comprehensive calculation of noise levels at the source ($\mathrm{NLS}(\cdot)$), transmission loss ($\mathrm{TL}(\cdot)$), and fuel consumption rate ($\mathrm{F}(\cdot)$) are explained in the subsequent sections.

\subsection{Acoustic model}
\label{sec_urn}

This section addresses the modeling URN from ships and introduces the first objective function. The noise level at the source (NLS) or near-field noise levels of the ship is the sound pressure level (SPL) measured at a reference distance of one meter from the source. It is independent of the environment in which the source operates. Meanwhile, the noise level received by mammals (NL), or far-field noise levels, is lower than the NLS due to transmission loss. It can be determined using a passive sonar equation, as follows:
\begin{equation}
\begin{aligned}
    \mathrm{NL} = \mathrm{NLS} - \mathrm{TL},
\end{aligned}
\label{eq_sonar}
\end{equation}
where all three quantities are in dB ref $1\mu$Pa. The transmission loss is the dissipation of energy influenced by several factors, including the loss attributed to geometric spreading, boundary effects, scattering phenomena, and volume attenuation. The propagation models compute acoustic TL for a fixed source location in the entire range-depth plane at all receiver depths and ranges.


Given that we know the received sound pressure level experienced by the mammals from Eq. \ref{eq_sonar}, we can introduce an objective function that quantifies the scope of the impact of URN. To encapsulate the impact of the entire voyage in a singular metric, the total acoustic intensity of the noise signal in units of W/m$^2$ is used as the objective function. The total sound intensity is given by
\begin{equation}
    \begin{aligned}
        I_{total} &= \int_{-\infty}^{\infty} I(f) \ df \\
        &= \sum_{f} I(f) \ df = \sum_{f} I_0 10^{\frac{SPL(f)}{10}} \ df.
    \end{aligned}
    \label{eq_intensity_power}
\end{equation}
The characteristic acoustic impedance of seawater is $1.5 \times 10^6$kg m$^{-2}$ s$^{-1}$, as a result, $I_0 = 0.67 \times 10^{-18}$ W/m$^2$ (i.e., 0 dB ref $1\mu$Pa) \cite{KUPERMAN2003317}.


Research indicates that an audiogram is a reliant reference, significantly above which can cause residual hearing effects on mammals \cite{Southall2019}. Following that, the intensity levels that exceed the specific mammal's audiogram threshold are determined by subtracting the SPL corresponding to the audiogram at each frequency. Subsequently, $I_{total}$ of the resulting SPL is calculated based on Eq. \ref{eq_intensity_power}. This quantifies the intensity of noise levels that surpass the threshold for an individual mammal at a given ship location. This summation is performed for all mammals and at each ship location, resulting in a unified metric that represents the ship's overall noise impact as written below
\begin{equation}
\begin{aligned}
\mathcal{J}_1 &= \sum_{i=1}^{n-1} \sum_{j=1}^{m} \sum_{f} I_0 10^{\frac{\mathrm{NL}_{ij}(f) \ - \ \alpha_{j}(f)}{10}} \ df .
\end{aligned}
\label{eq_obj1}
\end{equation}

To sum up, the noise levels from the ship on $i^{th}$ sailing leg received by the $j^{th}$ mammal, $\mathrm{NL}_{i,j}$, is calculated from Eq. \ref{eq_sonar}. Then, the total acoustic intensity of noise levels exceeding the audiogram threshold ($\alpha_{j}$) is calculated and summed across all mammals and ship locations based on Eq. \ref{eq_obj1}.
The subsequent sections explain the modeling of NLS and TL for a given environment in order to get the NL as well as the audiogram functions.

\subsubsection{Near-field noise levels}
\label{sec_nearF}
The near-field refers to the region in close proximity to the source (ship), where the amplitudes of the sound field are influenced by the physical dimensions of the sound source itself. Importantly, the near-field is independent of the surrounding environment. 

Ship noises are attributed to various factors, including the propeller parameters, onboard machinery, and the movement of the hull as it traverses through the water. While the characterization of hydro-acoustic noise generated by ships remains limited in understanding, it is evident that, in general, reducing the speed of a ship with a fixed pitch propeller leads to a reduction in overall noise levels. Over the years, numerous empirical models have been proposed with the objective of modeling the underwater noise spectrum emitted by a ship, considering factors such as its speed and other relevant characteristics \cite{Zhu2022, Wittekind2014}.

Ross's models \cite{Ross1976}, derived from measurements of ships equipped with fixed-pitch propellers, incorporate both the speed ($v$) and displacement ($dt$) of the vessel in their formulas. Various other sophisticated source-level models specific to particular ships exist \cite{chion2019meta} and can be used to replace Ross's model. Ross's model gives generalized expected NLS in dB ref $1\mu$Pa. It is formulated as follows:
\begin{equation}
\begin{aligned}
    \mathrm{NLS}(f;v,dt) =& 112 + 50 \log(\frac{v}{10}) + \\
    & 15\log(dt) + 20 - 20\log(f),
\end{aligned}
\label{eq_ross}
\end{equation}
where $v$ is ship speed (kt), $dt$ is displacement (MT), and $f$ is frequency (Hz). The above equation is used to calculate the noise levels 
($\mathrm{NLS}_i$) at each sailing leg based on the corresponding ship velocity ($v_i$) as shown in Fig. \ref{fig_net}.

\begin{figure}[ht]
\centering
\includegraphics[width=7.8cm]{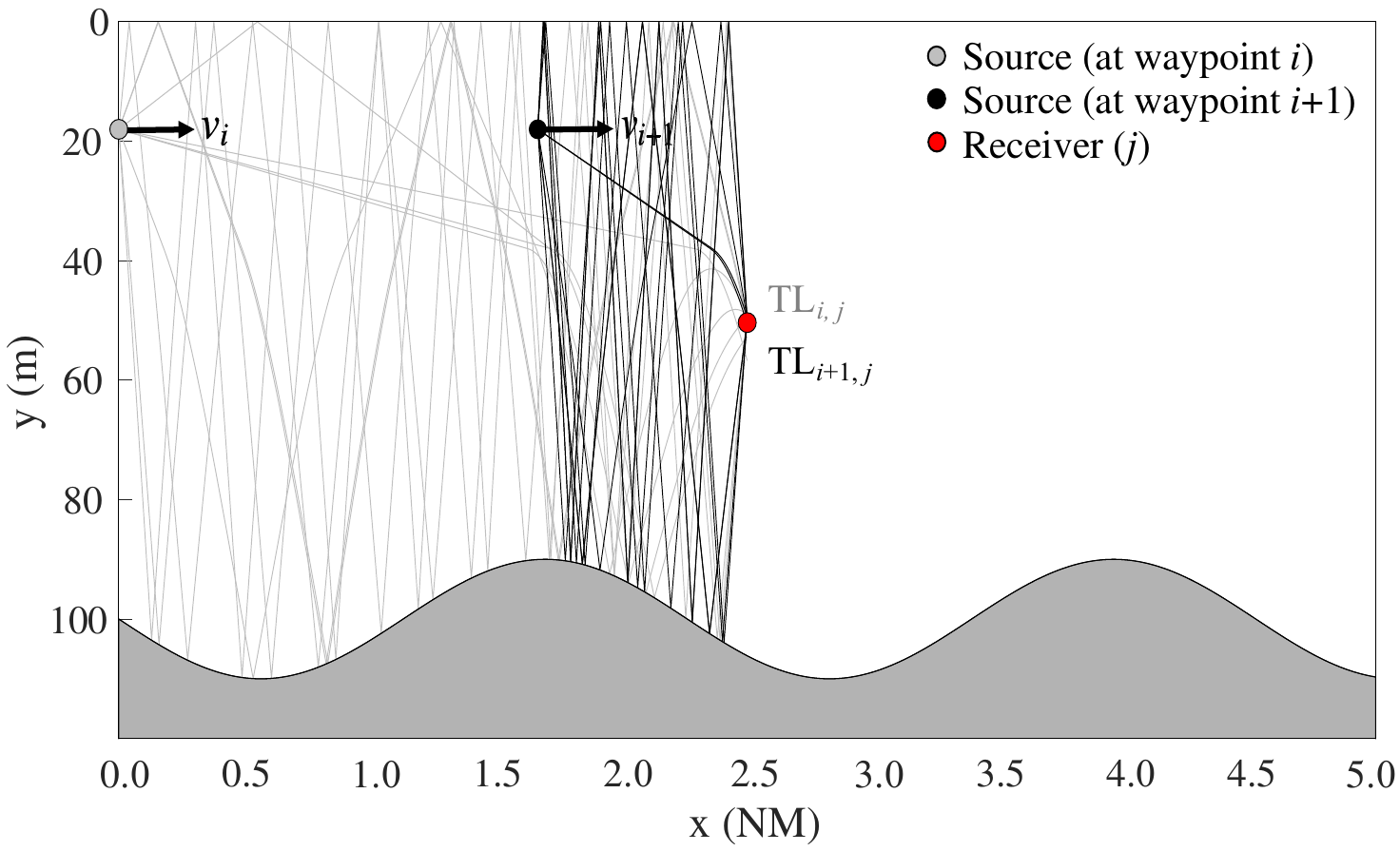}
\caption{Eigenrays plot from a moving source representing a ship to a fixed receiver representing a marine mammal. Only a few eigenrays are plotted for clarity. }
\label{fig_eigenplot}
\end{figure}

\subsubsection{Far-field noise levels}
\label{sec_farF}
The far-field is the region where the source can be approximated as a point source. The pressure amplitude monotonically decreases as the distance from the source increases. Moreover, these propagation losses are significantly influenced by the environment. Various models for acoustic propagation exist, based on different methods including ray tracing, normal modes, parabolic equations, and wavenumber integration. Recently, data-driven models generalized for varying complicated ocean environments have emerged, demonstrating accurate predictions of propagation characteristics \cite{mallik2024deep, Mallik2022}.

The far-field noise level in this study is obtained by using the Gaussian beam method. This method begins with the integration of the ray equations to obtain the central ray of the beam. Beams are then
constructed about the central ray by integrating a pair of auxiliary
equations, which govern the evolution of the beam in terms
of the beamwidth and curvature as a function of arc length. In a cylindrical coordinate system with r denoting the horizontal range and z the depth coordinate. The ray equation is given by:
\begin{equation}
\begin{aligned}
\frac{d}{d s}\left(\frac{1}{c(r, z)} \frac{d \mathbf{r}}{d s}\right)=-\frac{1}{c^2(r, z)} \nabla c(r, z),
\end{aligned}
\label{eq_ray}
\end{equation}
where $\mathbf{r}=\mathbf{r}(s)$, the $[r(s), z(s)]$ coordinate of the ray as a function of the arc length $s$, and $c(r, z)$ is the sound speed. This equation is reduced to the first-order by introducing the auxiliary variables $(\rho, \zeta)$. The detailed derivation of the Gaussian beam equations is given in the article by \cite{porter1987gaussian}. 

Now, in order to obtain the pressure field, each ray is assigned a phase and amplitude. The sound pressure amplitude along the ray is governed by a transport equation and is solved using the idea of the geometrical spreading of the ray tube. The variation in intensity along a ray tube is determined by its inverse relation with the cross-sectional area of the tube. On the other hand, the phase of the ray is computed based on the time of travel along the ray. For more details, the reader is referred to \cite{cerveny1987ray}.

At any given point, the pressure field is calculated by first identifying eigenrays that is, every ray that passes through that point from the source. The eigenrays from a source to a receiver are depicted in Fig. \ref{fig_eigenplot}. The pressure field is contributed by each of these eigenrays ($p_j(r, z)$) based on its amplitude and phase. The total sound pressure can be expressed as:
\begin{equation}
\begin{aligned}
p(r, z)=\sum_{j=1}^{N(r, z)} p_j(r, z),
\end{aligned}
\label{eq_pressure}
\end{equation}
where $N(r, z)$ is the number of eigenrays of the receiver. Hence the transmission loss computed in dB at the receiver is given by:
\begin{equation}
\begin{aligned}
\mathrm{TL}=-20 \log \left|\frac{p(r,z)}{p_0}\right|,
\end{aligned}
\label{eq_transmissionloss}
\end{equation}
where $p_0$ is the pressure at 1 m distance from the source. As shown in Fig. \ref{fig_eigenplot}, the source moves from one waypoint ($i$) to another ($i+1$). Consequently, the transmission loss (TL$_{i, j}$) at a fixed receiver ($j$) from the source at each waypoint ($i$) is computed based on the equation above. The computation is performed for frequencies from 10 to $10^4$ Hz. Furthermore, the entire process is repeated for all the receivers, which are the scattered marine mammals. The next section quantifies the impact of URN on these mammals.

\subsubsection{Marine mammal's audiogram functions}
\label{sec_mammal} 
 Numerous studies have been conducted to quantify and objectively evaluate the noise levels that may cause temporary or permanent threshold shifts in various marine mammals \cite{Weilgart2007, Southall2019}. One crucial criterion characterizing the hearing and absolute sensitivity of an auditory system is the audiogram. An audiogram represents a U-shaped curve of detection thresholds across a range of frequencies for a particular mammal or group. While audiograms come with extensive assumptions, extrapolation, and caveats, such as not accounting for variations in sensitivity with depths, gender, and age dependence, they still serve as a reasonable hearing demographic to start with. This criterion can be further refined by incorporating audiogram weighting functions, noise exposure levels, and other relevant metrics. As a result, the audiogram is regarded as the critical threshold, and noise levels surpassing this limit are interpreted as potentially harmful to marine mammals. 

Marine mammals are segregated into groups of similar hearing characteristics and frequency ranges. There are a total of seven hearing groups: High-Frequency Cetacean (HF), Very-High-Frequency Cetacean (VHF), Sirenian (SI), Phocid Carnivores in Air (PCA) and Water (PCW), Other Marine Carnivores in Air (OCA), and Water (OCW). Each of these groups is represented by a data-fitted parameterized mammalian audiogram function \cite{Southall2019}. Thus, the audiogram function of the $j^{th}$ mammal is formulated as follows,
\begin{equation}
\alpha_{j}(f)=\alpha_{0} + P_1 \log _{10}\left(1+\frac{P_2}{f}\right)+\left(\frac{f}{P_3}\right)^{P_4},
\label{eq_audiog}
\end{equation}
where $\alpha_{0}$,$P_1$, $P_2$, $P_3$, $P_4$ and $P_5$ are the parameter values of a hearing group. Since this study focuses specifically on the underwater acoustic propagation model, the consideration of marine carnivore groups in the air has been omitted. Consequently, we have a total of five distinct marine mammal hearing groups and corresponding functions. These functions represent each mammal's audiogram and are used in Eq. \ref{eq_obj1}.
The subsequent section presents the fuel consumption model and introduces the second objective function of the optimization problem.

\subsection{Fuel consumption model}
As mentioned before, the fuel consumption rate for the $i^{th}$ sailing leg is represented by $\mathrm{F}(v_i;\phi)$ in MT/h, where, $v_i$ and $\phi$ is the ship velocity and known set of ship characteristics, respectively. Therefore, the total fuel consumption of the voyage which is the second objective function is given by:
\begin{equation}
\begin{aligned}
    \mathcal{J}_2 = \sum_{i=1}^{n-1} \frac{d_i}{v_i} \mathrm{F}(v_i;\phi).
\end{aligned}
\label{eq_obj2}
\end{equation}

The fuel consumption rate can be expressed as follows:
\begin{equation}
\begin{aligned}
    \mathrm{F}(v_i;\phi) &= P_B(v_i;\phi) \Psi_{SFOC} , \\
    P_B(v_i;\phi) &=  \frac{P_E(v_i;\phi)}{\eta_D (v_i;\phi)}, \\
    P_E(v_i;\phi) &= R_{total}(v_i;\phi)  v_i,
\end{aligned}
\label{eq_power_sfoc}
\end{equation}
where $P_B$ is the required brake power of the ship (W), $\Psi_{SFOC}$ is the specific fuel oil consumption (MT/kWh), $P_E$ is the effective power (W), $\eta_D$ is the total propulsion efficiency, and $R_{total}$ is the total water resistance (N) experience by the ship which is calculated using Lap-Keller method \cite{Keller1973}. The SFOC curve which depicts the relationship between fuel consumption rate and engine load, varies across different ships depending on their respective engines. \ref{sec_props_ship_env} provides further details on the SFOC curve.

The Lap-Keller graphs are used to calculate the still water resistance, primarily influenced by factors such as the ship's speed, displacement, and hull form \cite{Keller1973}. This empirical approach takes into account various types of resistance, including frictional resistance, (caused by friction along the wetted surface), residual resistance (encompassing wave resistance, viscous pressure resistance, and resistance due to hull curvature), and incremental resistance (arising from surface roughness of the hull). The Lap-Keller method can be replaced with more robust and accurate models to estimate the ship's performance under dynamic sea conditions. Semi-empirical models \cite{lu2015semi} or data-driven models \cite{farag2020development} can used to account for additional resistance caused by the real state of the sea environment. Moreover, studies on ship voyage optimization have considered ocean currents \cite{yang2020ship}, metocean forecast models (including their uncertainties) \cite{wang2020effectiveness}, and even formulated these models as convex functions for faster real-time computation \cite{van2023convex}. These ideas can be incorporated into this general framework.

The total water resistance based on this method is expressed as follows:
\begin{equation}
\begin{aligned}
     R_{total} =  R_f + R_r + R_i = [C_f + C_r + C_a] \frac{1}{2} \rho v^2 S,
\end{aligned}
\end{equation}
where $C_f$, $C_r$, and $C_a$ are the specific frictional, residuary, and incremental resistance coefficients, respectively. $\rho$ is the  seawater density (kg$/\text{m}^3$), and $S$ is the wet surface area ($\text{m}^2$). \ref{sec_lap_verif} provides the validation of the Lap-Keller method. 
With the objective functions defined, the next section defines the optimization problem and the procedure used to solve it.

\section{Optimization framework}
\label{sec_opt_toolbox}
The multi-objective optimization framework is designed for a ship that follows a predetermined fixed route, aiming to minimize two objective functions defined in the previous sections:
i) total noise intensity levels ($\mathcal{J}_1$) and ii) total fuel consumption ($\mathcal{J}_2$). In this framework, the sailing speed of the ship is the control variable, which is optimized while considering the operational conditions and adhering to the voyage constraints. The operation conditions include the velocity limits and maximum engine load restriction while ensuring that the estimated time of arrival (ETA) is also satisfied.

Normalization of objective functions is necessary to make the objective comparable; otherwise, it leads to biased sampling on the Pareto front. The two objective functions from Eq. \ref{eq_obj1} and Eq. \ref{eq_obj2} are normalized using the corresponding ideal ($\mathcal{J}^*_i$) and nadir ($\mathcal{J}^N_i$) objective vectors which is defined below \cite{Kalyan2001}:
\begin{equation}
\begin{aligned}
\mathcal{J}^*_i &= \{\mathcal{J}_1(v^{[1]}), \mathcal{J}_2(v^{[2]})\} \\
&\text { where } \ \ \mathcal{J}_i(v^{[i]}))=\min_{v \in R^{n-1}}\left\{\mathcal{J}_i(x)\right\} \\
\mathcal{J}^N_i &= \{\mathcal{J}_1(v^{[2]}), \mathcal{J}_2(v^{[1]})\}.
\label{eq_ideal_nadir}
\end{aligned}
\end{equation}
Hence, the objective functions are normalized as follows:
\begin{equation}
    \hat{\mathcal{J}}_i = \frac{\mathcal{J}_i - \mathcal{J}^*_i}{\mathcal{J}^N_i - \mathcal{J}^*_i}.
    \label{eq_obj_norm}
\end{equation}

The optimization problem constituted by the aforementioned normalized objective functions and inequality constraints is written as
\begin{equation}
\begin{aligned}
\min_{v \in R^{n-1}} &  \{ \hat{\mathcal{J}}_1,\hat{\mathcal{J}}_2 \}, \\
\textrm{s.t.}
& \left\{\begin{array}{l}
\sum_{i=1}^{n-1} \frac{d_i}{v_i}-ETA \leq 0 \\
v_i-v_{\max } \leq 0, i=1, \cdots, n-1 \\
v_{\min }-v_i \leq 0, i=1, \cdots, n-1 \\
P_i-P_{\max } \leq 0, i=1, \cdots, n-1
\end{array}\right. \\
\end{aligned}
\end{equation}
where $\hat{\mathcal{J}}_1$ and $\hat{\mathcal{J}}_2$ are the two normalized objective functions defined as in Eq. \ref{eq_obj_norm}, $v_i$ is the velocity of the ship at $i^{th}$ sailing leg (kt), $v_{\min}$ and $v_{\max}$ is the minimum and maximum ship velocity (kt), respectively, $P_i$ is the total engine power at $i^{th}$ sailing leg (kW), and $P_{\max}$ is the maximum engine power of the ship (kW).

The genetic algorithm is a population-based meta-heuristic algorithm commonly employed to solve optimization problems \cite{Beheshti2013}. The non-dominated sorting genetic algorithm (NSGA-II) is based on the genetic algorithm and is utilized for solving MOOPs \cite{Verma2021}. NSGA-II generates an evenly distributed Pareto front with robust convergence and is widely acclaimed in the field of MOOPs. NSGA-III is a further advancement over NSGA-II for handling many objectives, introducing a new reference-point-based selection strategy that enhances population diversity. It has been proven effective for test problems with three or more objectives \cite{deb2013evolutionary}; however, our framework has only two objective functions. Hence, NSGA-II is applied to address the aforementioned MOOP, generating a collection of optimal ship speeds collectively referred to as the Pareto front.

NSGA-II starts with randomly sampling a chromosome, which constitutes a population set. Then the chromosomes are evaluated using the objective function, followed by three operations: (i) selection, (ii) crossover, and (iii) mutation. In this paper, binary tournament selection is implemented, where two chromosomes from the population set are randomly selected and compared head-to-head to be added to the parent set. It is the most popular selection method, offering advantages such as time efficiency and diversity preservation \cite{goldberg1991comparative}. Next, the crossover operation is applied to the parent set to produce better offspring for the next generation. A simulated binary crossover is used in this paper as it is better suited for MOOPs \cite{umbarkar2015crossover}. Furthermore, the mutation operation is applied, which randomly picks a chromosome and alters its genes. The polynomial mutation is used to enhance the exploration of the search space.

In this study, the NSGA-II algorithm is implemented using the pymoo library in Python \cite{pymoo}, which is a framework for both single and multi-objective optimization, offering state-of-the-art optimization algorithms. The parameters of NSGA-II include population size, total number of generations, crossover probability, and mutation probability. The appropriate setting of these parameters for evolutionary algorithms is crucial for achieving best performance. The task of parameter tuning is both tedious and time-consuming, necessitating the use of various automatic parameter-tuning approaches \cite{huang2019survey}. Among these, the Sequential Model-Based Algorithm Configuration (SMAC) stands out as one of the most powerful methods for parameter tuning. We employ SMAC3 for the parameter tuning of NSGA-II \cite{smac3}, utilizing this open-source Python package that incorporates a Bayesian optimization approach combined with a random forest or Gaussian process model to estimate the performance of different configurations. The hypervolume indicator is used as the performance metric in this algorithm to select optimal parameters \cite{shang2020survey} suitable for this application. Further details are provided in \ref{sec_hyperparam}.






Once the Pareto optimal set has been acquired, multiple criteria decision-making (MCDM) methods are employed to systematically assess and identify a trade-off solution among the Pareto optimal solutions \cite{Sahoo_Goswami_2023}. The Technique for Order of Preference by Similarity to Ideal Solution (TOPSIS) is a simple and widely used classical MCDM method that determines trade-off solutions by ranking all the solutions on the Pareto front \cite{hwang2012multiple}. It ranks solutions based on their proximity to the positive ideal solution (PIS) and their distance from the negative ideal solution (NIS). In this study, the PIS is the solution with the lowest noise intensity levels and fuel consumption, whereas the NIS represents the maximum levels of these two attributes.

Initially, a decision matrix is constructed by stacking the solutions of the Pareto front row-wise, with the corresponding attributes (noise intensity levels and fuel consumption) positioned column-wise. Subsequently, we normalize the decision matrix using a vector normalization technique and assign equal weights to both attributes. Many normalization techniques exist; however, this specific one was chosen because it has proven to be effective and suitable for the TOPSIS method \cite{vafaei2018data}. Then, the PIS and NIS are determined from the normalized decision matrix, followed by the measurement of separation for each solution from these two points. Finally, the relative closeness coefficient for each solution is calculated based on the separation and ranked in descending order. The highest-ranked solution is considered the trade-off solution.



\begin{figure}[h!]
\centering
\includegraphics[width=7.8cm]{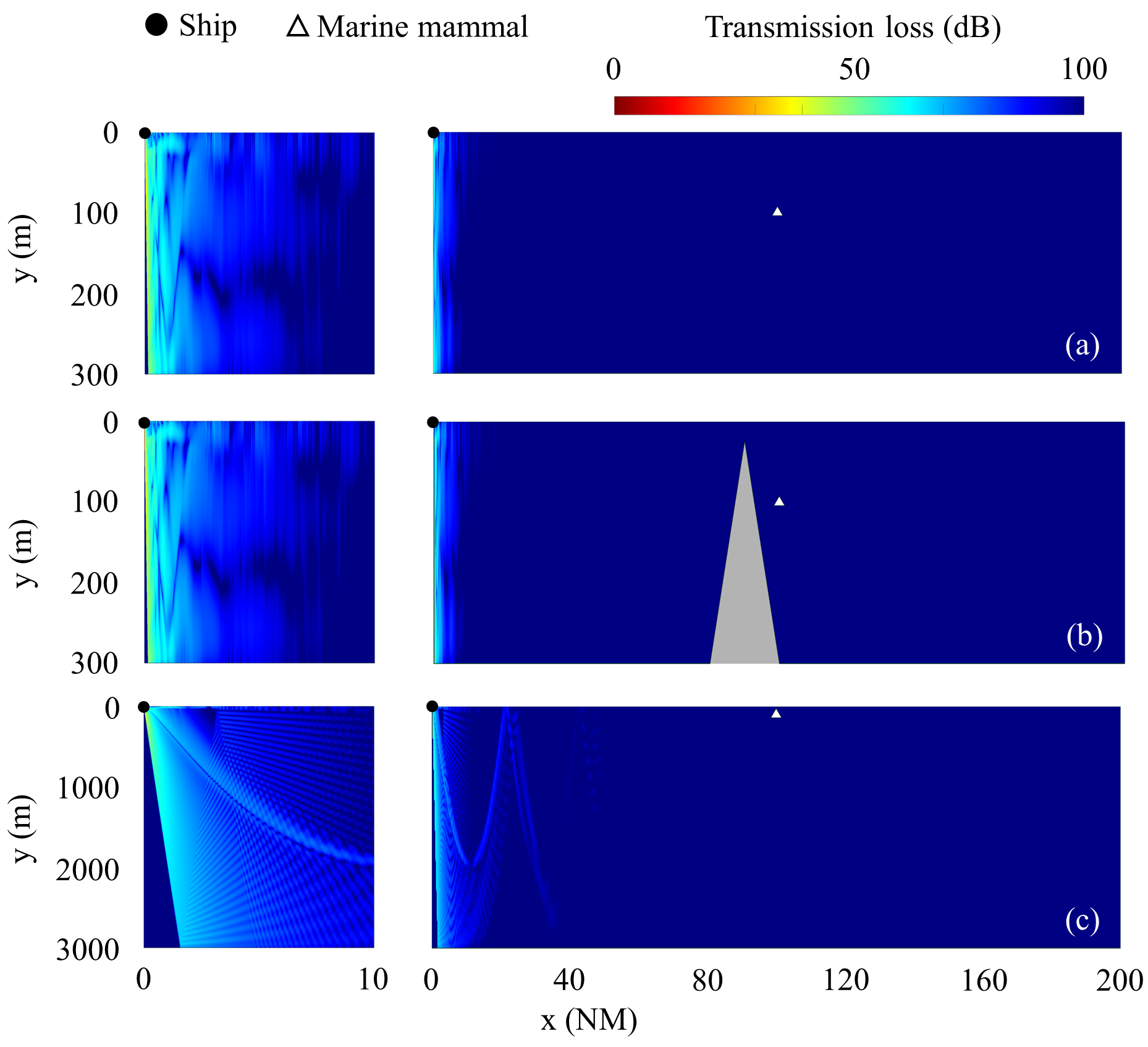}
\caption{Configurations for sub-cases in Case T: (a) Case T1, (b) Case T2, and (c) Case T3. }
\label{fig_case_T}
\end{figure}

\section{Case studies} \label{sec_case_studies}

In this section, we present the case studies to illustrate the applicability of the framework. To provide a comprehensive analysis, we present three primary case studies: Case T (Test), Case A (Shallow water), and Case B (Deep water).

Case T (Test) highlights the impact of bathymetry and the presence of single marine mammals precisely positioned at the voyage's center. Consequently, the following sub-cases are examined:

\begin{itemize}
    \item T1: Incorporation of a single mammal at the center of the ship's path in the shallow-water region.
    \item T2: Incorporation of a single mammal at the center of the ship's path with distinctive bathymetry (seamount), to provide elucidation of bathymetry's effects.
    \item T3: Incorporation of a single mammal at the center of the ship's path in the deep-water region.
\end{itemize}

 \begin{figure*}[ht]
\centering
\includegraphics[width=14cm]{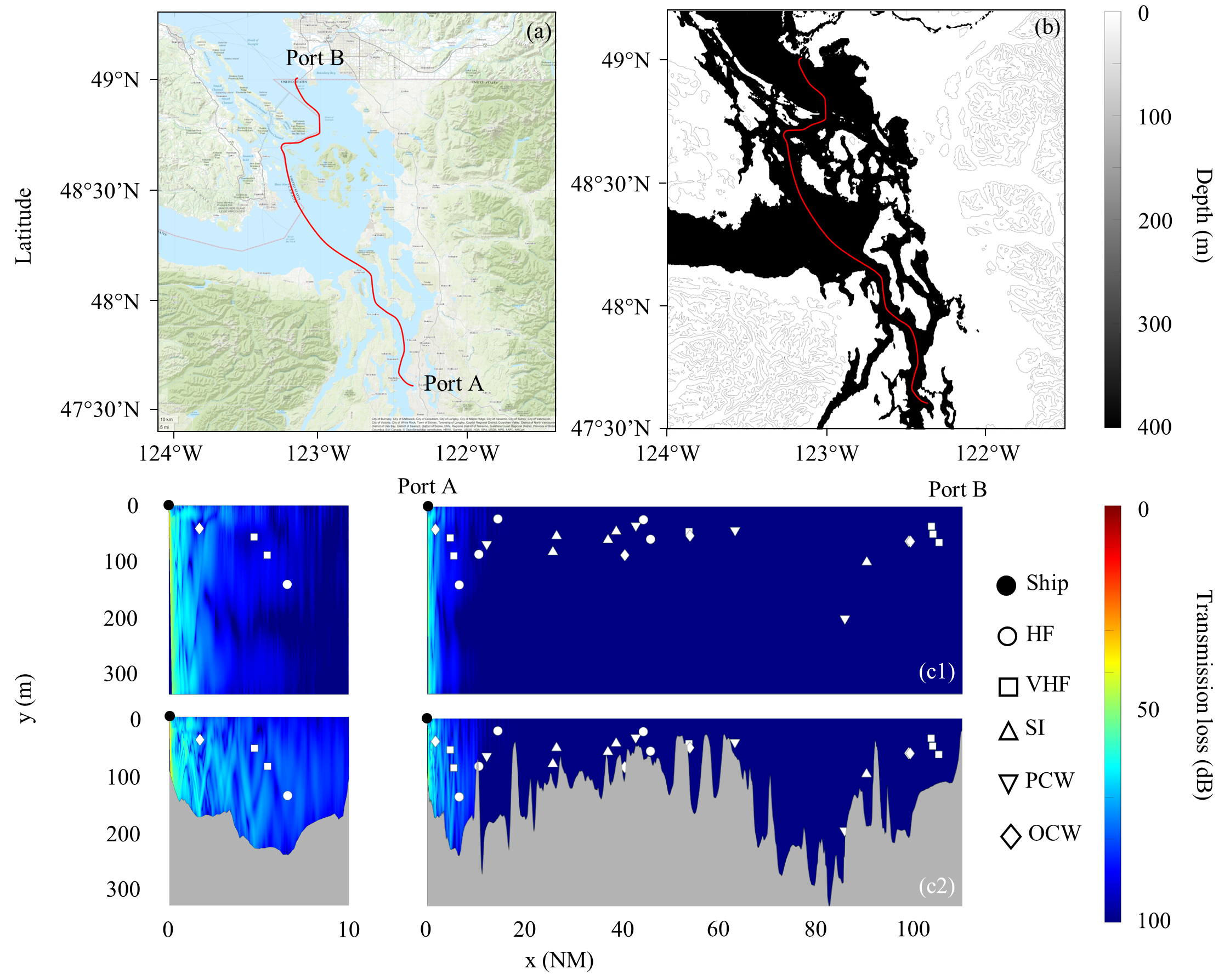}
\caption{Configurations for sub-cases in Case A: (a) Ship route on the geographic plot; (b) Bathymetry contour; (c) Transmission loss contour along with a 10 NM zoom box for (c1) Case A1, and (c2) Case A2. }
\label{fig_Case_A}
\end{figure*}

On the other hand, Case A (Shallow water), and Case B (Deep water), highlight the differentiation between shallow-water and deep-water regions. In the shallow environment, sound encounters multiple reflections and scattering from both the sea floor and the surface, in contrast to the deep environment. Consequently, shallow-water sound propagation is heavily influenced by the structure of the sea floor. Additionally, considering the extensive presence of continental shelves in vast oceanic areas, characterized by shallow water, it becomes imperative to account for such regions. As a result, our case study is demarcated into two distinct regions of interest,i.e., Case A (Shallow water), and Case B (Deep water).


\begin{figure*}[ht]
\centering
\includegraphics[width=14cm]{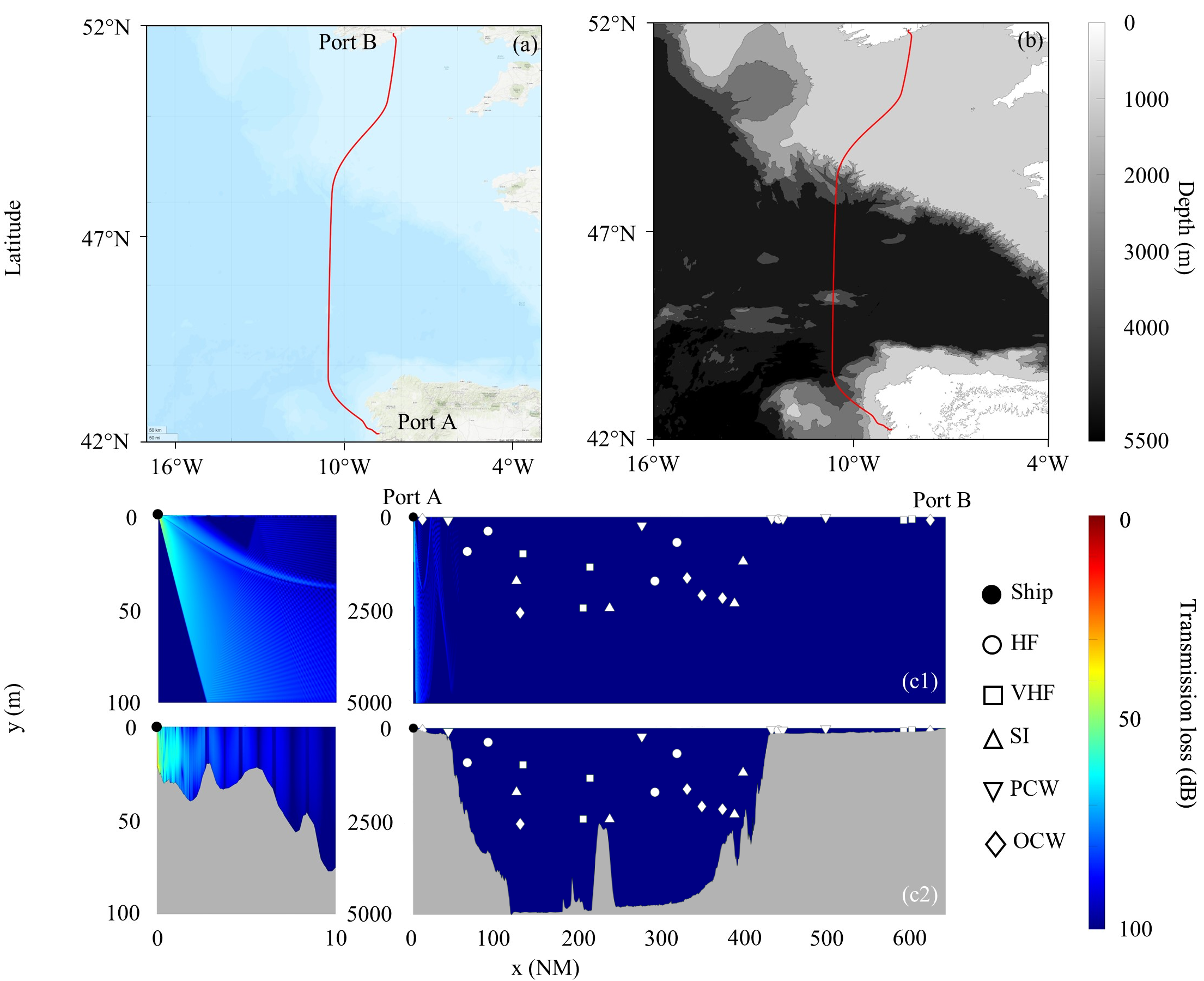}
\caption{Configurations for sub-cases in Case B: (a) Ship route on the geographic plot; (b) Bathymetry contour; (c) Transmission loss contour along with a 10 NM zoom box for (c1) Case B1, and (c2) Case B2. }
\label{fig_Case_B}
\end{figure*}

Case A presents the simulation of a 6900 TEU containership traveling from the Port of Seattle, United States (Port A) to Roberts Bank Port, Canada (Port B) spanning a distance of 110.09 NM. This route is chosen specifically because the Salish Sea is a shallow region. The ship's route spans a maximum depth of approximately 331.2 m and an average depth of 158.1 m. On the other hand, in Case B for the deep-water region, the same ship starts on a journey from Port of Vigo, Spain (Port A) to Port of Cork, Ireland (Port B), covering a total distance of 641.68 NM. The maximum and average depth along the path is around 5042.5 m and 503.1 m.

Case A and B are further divided as follows to gain a deeper understanding of the effects of real bathymetry data:
\begin{itemize}
    \item A1 \& B1: Incorporation of 25 marine mammals, randomly selected from five distinct hearing groups mentioned in Section \ref{sec_mammal} and randomly scattered along the path. Bathymetry is maintained at a constant depth along the entire route.
    \item A2 \& B2: Similar to Case A1 and B1, however, here the bathymetry data extracted at a resolution of 0.31 NM from the GEBCO \cite{gebco2023} database is used.
\end{itemize}

In order to proceed with the aforementioned cases, it is essential to clearly define the ship's operational and environmental conditions. The ship's propulsion power is dependent on the water resistance, as indicated in Eq. \ref{eq_power_sfoc}, which in turn is influenced by various parameters of the ship. All these necessary ship characteristics required for the calculation are in \ref{sec_props_ship_env} along with the environmental conditions. A pragmatic ship path from Port A to Port B is made to facilitate our analysis. Fig. \ref{fig_Case_A} and Fig. \ref{fig_Case_B} show the predetermined route in conjunction with the bathymetry contour of Case A and B, respectively. Note that the path and sound speed profile (SSP) remain fixed across all sub-cases. The ETA is set to 16, 10, and 60 hours for Case T, A, and B, respectively. 

 The initial states of all three primary cases are depicted in Fig. \ref{fig_case_T}, Fig. \ref{fig_Case_A}, and Fig. \ref{fig_Case_B}. The test case, Case T, includes a ship covering a total distance of 200 NM with a mammal at the center. The plot illustrates the TL at 20 Hz for three sub-cases along with a zoom box covering a 10 NM range. Whereas for Cases A and B, the ship's route on the geographic plot and bathymetry contour are presented side by side. Furthermore, below these figures, the TL at 20 Hz along with multiple mammals of varying audiogram groups are plotted for both sub-cases.

 It is evident in all these TL contours that the noise levels attenuate quickly over shorter distances, especially in the presence of bathymetry. However, a contrasting characteristic arises between Case A and Case B. Within the zoom boxes of Fig. \ref{fig_Case_A} (c2) and Fig. \ref{fig_Case_B} (c2), the noise levels of Case B2 fade quicker than Case A2. This is attributed to the fact that, in the case of shallow-water regions, noise signals reflect frequently, preventing them from penetrating deeper depths and allowing them to persist over longer distances than in deep-water regions. Furthermore, in Case B1, the refraction of sound waves can be recognized, resulting in shadow regions close to the surface. These observations underscore the complexity of noise signals propagating in the ocean, highlighting the necessity of adapting the ship to these environmental conditions. In the subsequent section,  the results of these case studies are discussed and compared.

\section{Results and discussion} \label{sec_results_d}
\begin{table*}[ht]
\centering
\begin{tabular}{ll|ccc|cc|cc}
\hline
\multicolumn{1}{c}{}                                                                                                               & \multicolumn{1}{c|}{} & \multicolumn{3}{c|}{Case T} & \multicolumn{2}{c|}{Case A} & \multicolumn{2}{c}{Case B} \\ \hline
\multicolumn{1}{c}{}                                                                                                               & \multicolumn{1}{c|}{} & T1      & T2      & T3      & A1           & A2           & B1           & B2          \\ \hline
\multicolumn{1}{l|}{\multirow{3}{*}{\begin{tabular}[c]{@{}l@{}}Noise dominant \\ solution ($\hat{\mathcal{J}}^*_1$)\end{tabular}}} & ETA (h)               & 15.63   & 15.39   & 14.77   & 9.99         & 9.98         & 49.39        & 50.29       \\
\multicolumn{1}{l|}{}                                                                                                              & $\mathcal{J}_1$ (dB)  & 6.53    & 6.53    & 6.14    & 9.75         & 9.80         & 7.95         & 8.20        \\
\multicolumn{1}{l|}{}                                                                                                              & $\mathcal{J}_2$ (MT)  & 2638.38 & 2648.46 & 2587.77 & 817.25       & 820.03       & 23564.57     & 23638.44    \\ \hline
\multicolumn{1}{l|}{\multirow{3}{*}{\begin{tabular}[c]{@{}l@{}}Fuel dominant\\ solution ($\hat{\mathcal{J}}^*_2$)\end{tabular}}}   & ETA (h)               & 16.00   & 16.00   & 16.00   & 10.00        & 10.00        & 60.00        & 60.00       \\
\multicolumn{1}{l|}{}                                                                                                              & $\mathcal{J}_1$ (dB)  & 8.52    & 8.67    & 7.12    & 11.34        & 11.26        & 9.05         & 9.18        \\
\multicolumn{1}{l|}{}                                                                                                              & $\mathcal{J}_2$ (MT)  & 1640.89 & 1639.02 & 1642.10 & 692.27       & 693.79       & 11648.59     & 11554.32    \\ \hline
\multicolumn{1}{l|}{\multirow{5}{*}{\begin{tabular}[c]{@{}l@{}}Optimal\\ solution ($\hat{\mathcal{J}}_t$)\end{tabular}}}           & ETA (h)               & 15.98   & 15.99   & 15.99   & 9.99         & 9.99         & 59.99        & 59.99       \\
\multicolumn{1}{l|}{}                                                                                                              & $\mathcal{J}_1$ (dB)  & 7.17    & 7.12    & 6.13    & 9.96         & 9.95         & 7.97         & 8.23        \\
\multicolumn{1}{l|}{}                                                                                                              & $\mathcal{J}_2$ (MT)  & 1687.54 & 1685.22 & 1667.16 & 698.62       & 698.12       & 11693.19     & 11645.84    \\
\multicolumn{1}{l|}{}                                                                                                              & $\hat{\mathcal{J}}_1$ & 0.046   & 0.045     & 0.028   & 0.050        & 0.034        & 0.003        & 0.007       \\
\multicolumn{1}{l|}{}                                                                                                              & $\hat{\mathcal{J}}_2$ & 0.034   & 0.021     & 0.019   & 0.016        & 0.015        & 0.004        & 0.007       \\ \hline
\end{tabular}
\caption{Results of all the cases in terms of ETA and objective functions.}
\label{tab_results}
\centering
\end{table*}

\begin{figure*}[ht]
\centering
\includegraphics[width=15cm]{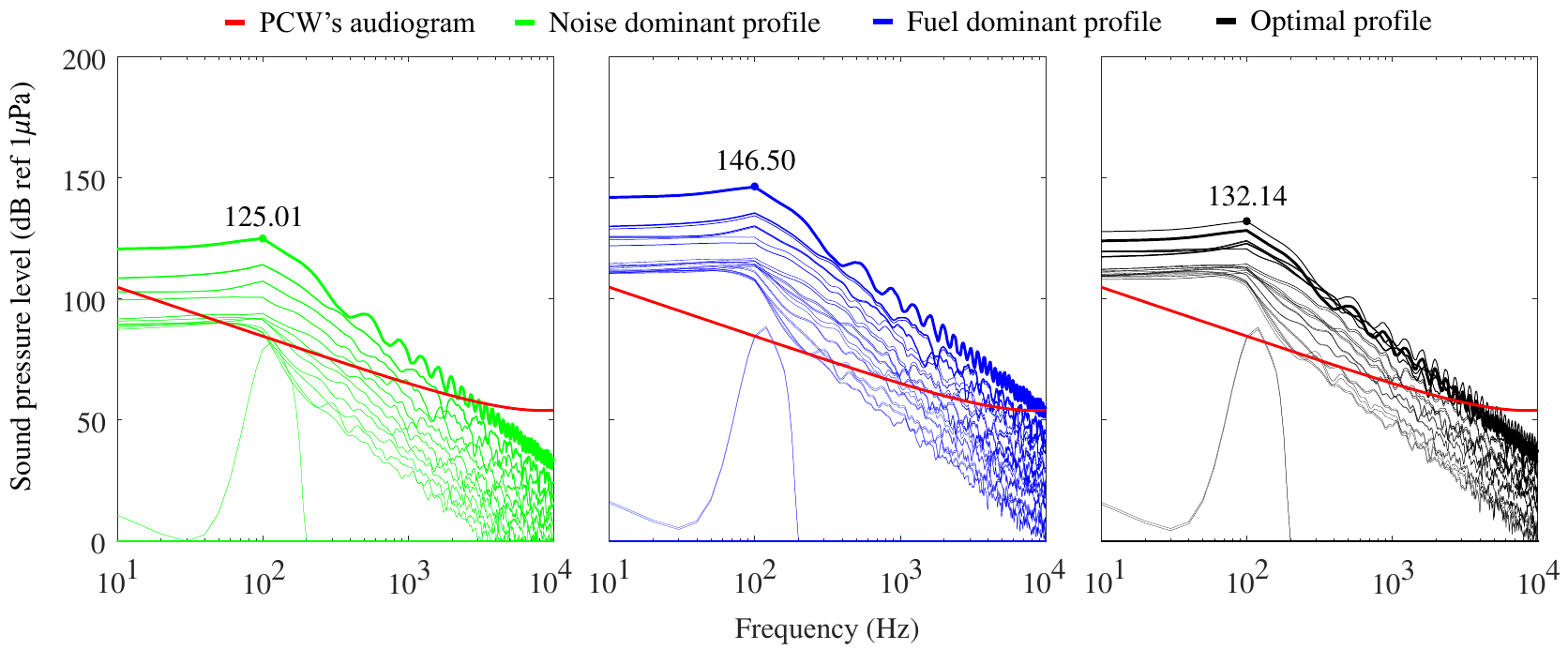}
\caption{Sound pressure levels received by the mammal from the ship on each sailing leg in Case T1, with thicker lines indicating proximity to the mammal.}
\label{fig_audiogram}
\end{figure*}

\begin{figure*}[ht]
\centering
\includegraphics[width=12.5cm]{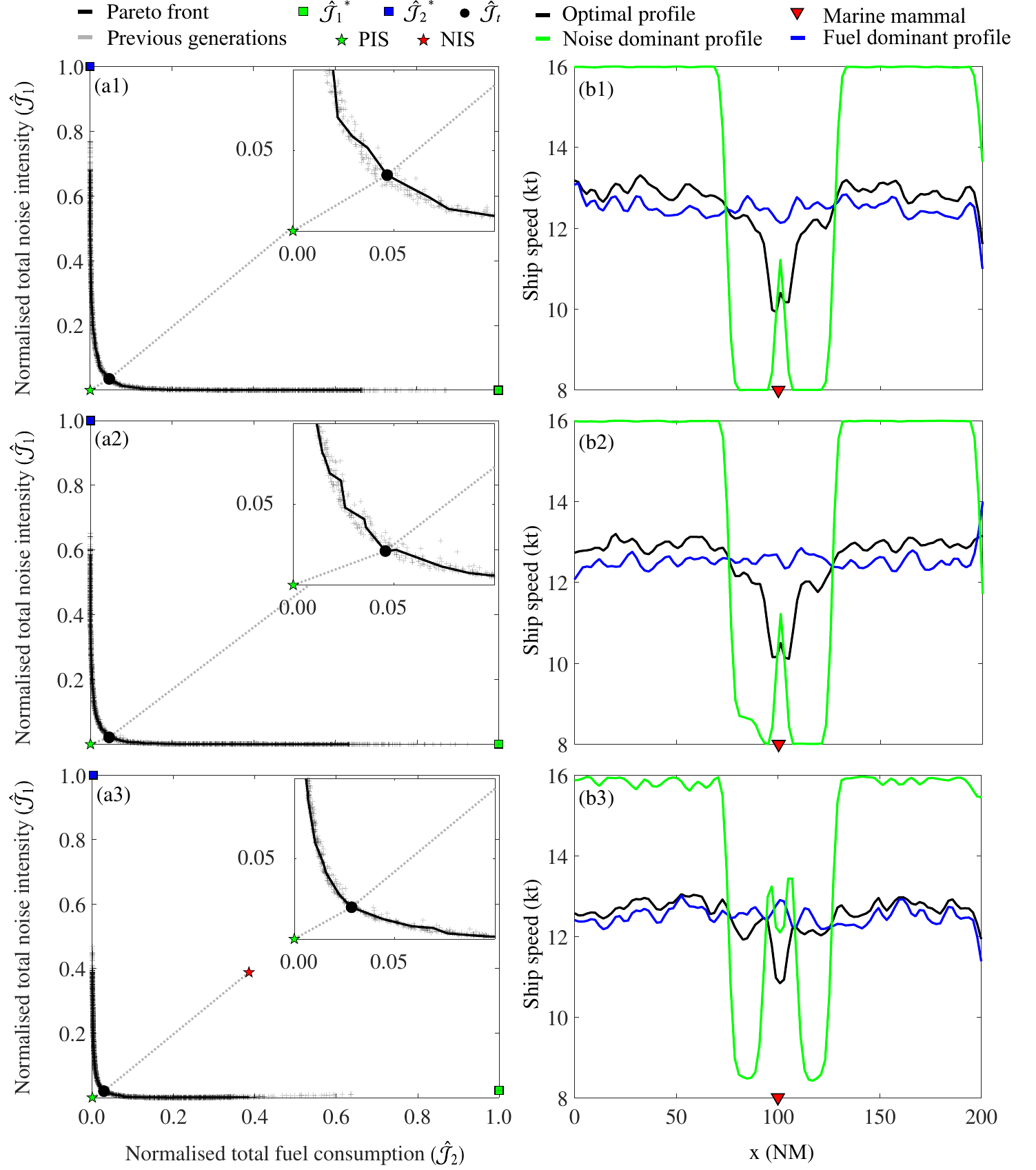}
\caption{Results of Case T: (a) Pareto front of $\hat{\mathcal{J}}_1$ vs $\hat{\mathcal{J}}_2$, and (b) Ship's sailing speed profile; (1) Case T1, (2) Case T2, and (3) Case T3.}
\label{fig_results_caseT}
\end{figure*}

\begin{figure*}[ht]
\centering
\includegraphics[width=14cm]{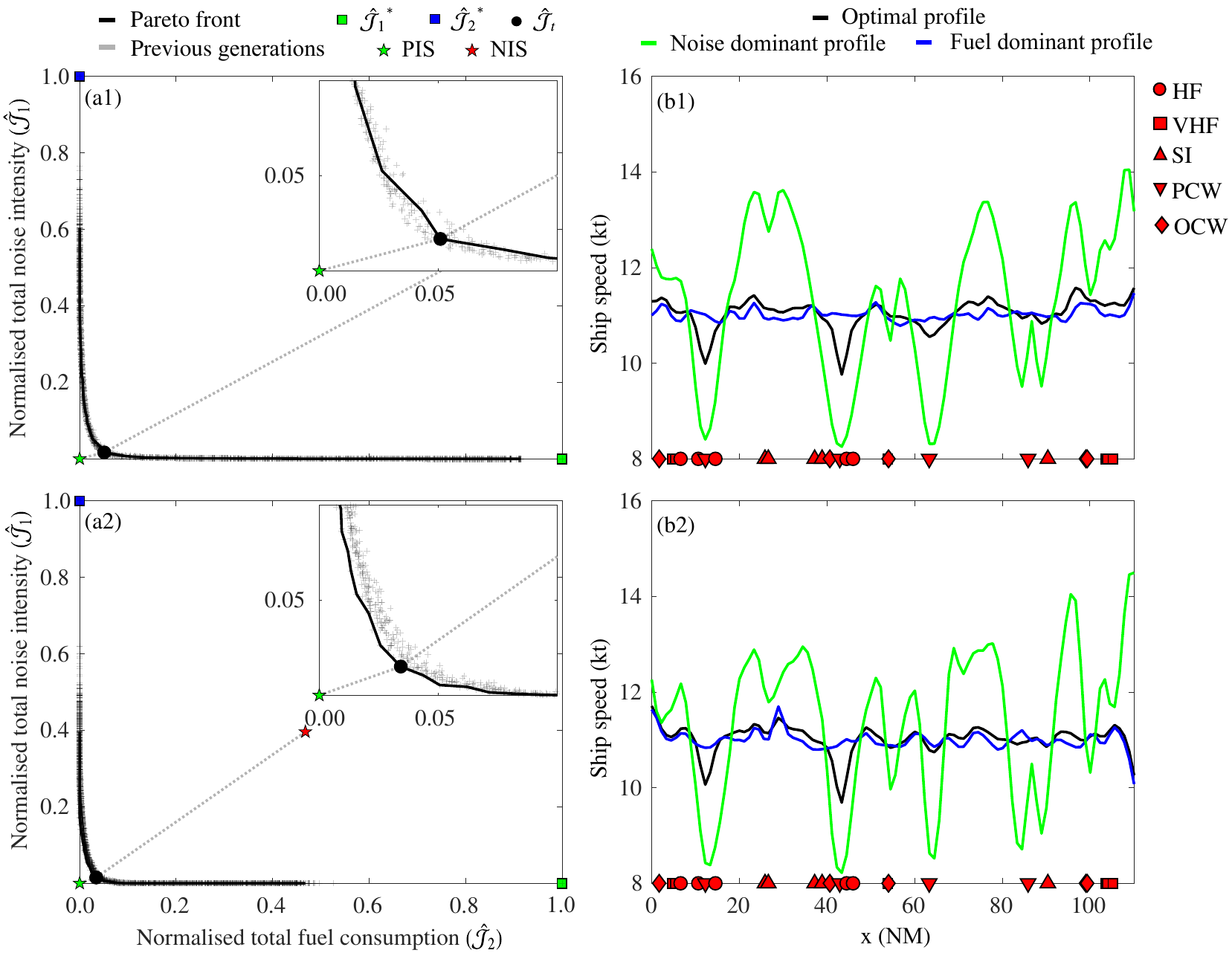}
\caption{Results of Case A: (a) Pareto front of $\hat{\mathcal{J}}_1$ vs $\hat{\mathcal{J}}_2$, and (b) Ship's sailing speed profile; (1) Case A1, and (2) Case A2.}
\label{fig_results_caseA}
\end{figure*}

\begin{figure*}[ht]
\centering
\includegraphics[width=14cm]{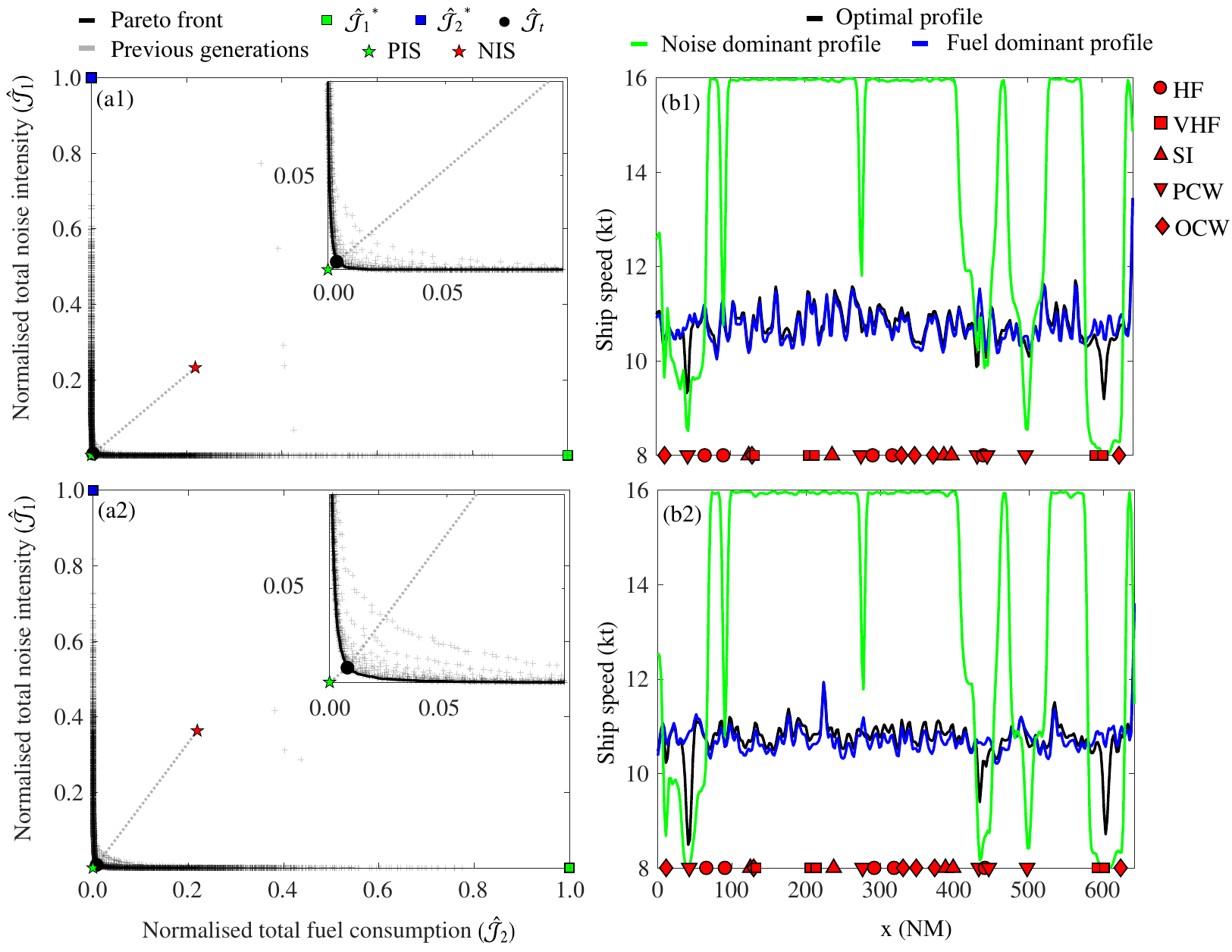}
\caption{Results of Case B: (a) Pareto front of $\hat{\mathcal{J}}_1$ vs $\hat{\mathcal{J}}_2$, and (b) Ship's sailing speed profile; (1) Case B1, and (2) Case B2.}
\label{fig_results_caseB}
\end{figure*}

The results comprise the normalized Pareto front and the optimized sailing speed profiles of the ship, presented side-by-side. The Pareto frontier encompasses multiple non-dominated solutions with five particular solutions being highlighted. The PIS and NIS from the TOPSIS method are indicated by a green star and a red star, respectively. The fuel-dominant ($\hat{\mathcal{J}}^*_2$) and noise-dominant objective vectors ($\hat{\mathcal{J}}^*_1$) (discussed in Section \ref{sec_opt_toolbox}) represented by the blue and green markers positioned at the extremes of the y- and x-axes, respectively.  Finally, the trade-off objective vector ($\hat{\mathcal{J}}_t$) between the two objective functions, estimated using TOPSIS, is represented by the black circle on the frontier as shown in the zoomed box. The solution vectors of every 200$^{th}$ generation are overlaid to show the convergence of the Pareto front. Moreover, the optimal trade-off speed profile of the ship is plotted next to the Pareto front in conjunction with the profiles corresponding to these highlighted objective vectors.

Before we embark on this, Fig. \ref{fig_audiogram} illustrates a series of sound pressure levels overlaid across the entire voyage of Case T1 for three distinct solutions: the noise-dominant, fuel-dominant, and optimal solutions. As mentioned before, a single mammal of the PCW type is positioned at the center of the path in Case T1. At each waypoint, the pressure level received by this mammal is computed according to the methodology outlined in Section \ref{sec_urn}. 
As expected, the noise-dominant solution has the lowest peak at 125.01 dB, in contrast to the fuel-dominant solution, which shows the highest peak at 146.50 dB. The optimal trade-off solution, on the other hand, shows its peak at 132.14 dB. Although the noise-dominant solution, optimized solely for the noise objective function, has the least sound pressure level compared to other solutions, it still surpasses the audiogram of the mammal at lower frequencies. This is due to exceptionally low audiogram levels, making it impractical to reduce below them under the voyage constraints. Next, the discussion focuses on the speed profiles and Pareto front of Case T.

Figure \ref{fig_results_caseT} shows the results of Case T. Here, the common trends on the normalized Pareto front of all the sub-cases are exponential convergence of the front along the x- and y-axis. The trade-off objective vector ($\hat{\mathcal{J}}_t$) in all three sub-cases falls within a margin of 0.05 for both normalized objective functions implying that the algorithm minimizes the two objectives by 95\% without compromising the other. 

Common trends in the speed profiles of Case T in Fig. \ref{fig_results_caseT} include the slowdown of the ship in the vicinity of the mammal which is around 27.28 NM away. As expected, the noise-dominant profile has a sudden decrease in speed close to the mammal, prioritizing to minimize noise. In contrast, the fuel-dominant profile exhibits no such slowdowns. The optimal profile shows a balanced compromise between these two characteristics. The optimal profile exhibits a higher voyage speed than the fuel-dominant profile both before and after the slowdown due to the time constraint. One of the peculiar features of the speed profile is the spike when the ship is exactly above the mammal. This is because of maintaining the angular range of the rays between -45 \textdegree to 45 \textdegree. Therefore, when the ship is precisely positioned above the mammal, minimal rays pass through the mammal due to sound refraction, leading to reduced noise levels and consequently, the spike.

Examining the sub-cases of Case T in Fig. \ref{fig_results_caseT} reveals speed profile oscillations in Case T3, which pertain to the deep-water region. As shown in Fig. \ref{fig_case_T}, deep water permits ray refraction, creating noise shadow and convergence zones. In deep waters, the upward refraction of sound leads to the spatially periodic formation of convergence zones with high intensity and shadow zones with low intensity, alternating near the surface. Therefore, the peaks of these oscillations are in correspondence with the locations of the shadow zones. As for Case T2, there is a kink in the noise-dominant profile due to the presence of the seamount immediately next to the mammal. Overall, the test cases outlined in Case T suggest that the optimized speed profile is influenced by the position of the marine mammal and bathymetry. 

A quantitative comparison of the optimized cost functions for all the sub-cases is tabulated in Table \ref{tab_results}. Focusing on the fuel-dominant solution ($\hat{\mathcal{J}}^*_2$) of Case T, the ETA (h) is satisfied across all three sub-cases, and total fuel consumption ($\mathcal{J}_2$ (MT)) remains consistent around 1640 MT. However, the noise objective function ($\mathcal{J}_1$ (dB)) exhibits the lowest value at 7.12 dB in Case T3. It should be noted that the fuel-dominant solution primarily optimizes $\mathcal{J}_2$, disregarding $\mathcal{J}_1$. Therefore, the lower value in Case T3 is attributed to the formation of shadow zones in the presence of deep bathymetry, as illustrated in Fig. \ref{fig_case_T}. In the comparison between Case T2 and T1, it is evident that the former exhibits higher noise levels, indicating a detrimental impact on the mammal positioned next to a seamount. As illustrated in Fig. \ref{fig_case_T}, while the seamount initially acts as a shield, as the ship traverses the peak of the seamount, it changes into a reflector, directing noise signals towards the mammal. This collectively results in a negative impact.

 Focusing on the noise-dominant solution ($\hat{\mathcal{J}}^*_1$) of Case T in Table \ref{tab_results}, again the ETA (h) is satisfied in all the sub-cases. Here, $\mathcal{J}_1$ is primarily minimized, neglecting $\mathcal{J}_2$. As expected, Case T3 exhibits the lowest noise levels at 6.14 dB due to shadowing effects. Despite achieving a reduction of 2 dB in $\mathcal{J}_1$ compared to the fuel-dominant solution, there is a significant increase of approximately 1000 MT in $\mathcal{J}_2$ compared to the fuel-dominant solution. This is the reason why we estimate the trade-off solution mentioned in Section \ref{sec_opt_toolbox} which prevents disproportionate fuel consumption for marginal reductions in noise power levels. The optimal solution ($\hat{\mathcal{J}}_t$) for Case T in the table indicates an approximate reduction of 1.5 dB in noise levels, followed by a slight increase in fuel consumption by 40 MT. In summary, the test cases highlight that a notable decrease in the overall noise intensity level can be achieved with only a minor increase in fuel consumption.

Moving on to Case A (Shallow water), the results are shown in Fig. \ref{fig_results_caseA}. The normalized Pareto front is similar to Case T where the optimal objective vector is within the 95\% range. The sailing speed profiles have some interesting features. Before discussing the speed profiles,  it is important to recognize that in this case marine mammals belonging to randomized audiogram groups are randomly scattered. As a result, clusters of mammals are found in close proximity in certain areas, while individual mammals are dispersed in other areas.

The optimized sailing speed for Case A is presented in Fig. \ref{fig_results_caseA}. In particular, the noise-dominant profile shows significant speed drops to 8 kt in the vicinity of mammal clusters, with minor decreases observed when encountering individual mammals. Furthermore, the decrease in speed is predominantly observed when the ship is near PCW-type mammals, implying that the optimization algorithm provides greater significance to this particular type. This preference is because PCW's audiogram has the lowest sound pressure level at low frequencies compared to other groups \cite{Southall2019}. $\mathcal{J}_1$ from Eq. \ref{eq_obj1} is explicitly formulated to account for intensity levels surpassing the audiogram. Consequently, this characteristic is inherited in the speed profile. In contrast, the fuel-dominant profile maintains a relatively constant speed, averaging around 11 kt with minor fluctuations. The optimal profile shares similarities with the fuel-dominant profile with some prominent characteristics of the noise-dominant profile. Comparing the noise-dominant speed profiles of Case A1 and A2 in Fig. \ref{fig_results_caseA}, it is evident that Case A2 exhibits more sharp and prominent features in the speed profile. This is attributed to bathymetry restricting prolonged sound propagation. Consequently, the optimization process becomes more localized in Case A2, in contrast to A1, which shows a globally smoother profile.

A quantitative comparison of Case A1 and A2 in Table \ref{tab_results} shows minimal differences. The noise objective function ($\mathcal{J}_1$) slightly favors A2, suggesting that the bathymetry acts as a protective barrier for marine mammals. Case A2, a realistic shallow-water case, demonstrates a reduction of 1.31 dB in total noise intensity (95.0\% lower compared to the fuel-dominant solution) through the optimization of sailing speed, requiring an additional 4.33 MT of fuel (0.6\% higher compared to the fuel-dominant solution).

Figure \ref{fig_results_caseB} shows the results of Case B (Deep water). The Pareto front demonstrates convergence similar to previous cases. One difference in this case is that the population size of each generation in NSGA-II was increased due to the voyage having more sailing legs, resulting in a greater number of decision variables, all to achieve convergence. The bathymetry of this case is plotted in Fig. \ref{fig_Case_B}, indicating shallow-water regions at the beginning and end; however, the rest of the voyage is in a deep-water environment. The noise-dominant speed profiles in Fig. \ref{fig_results_caseB} show patterns in conjunction with this bathymetry. The framework decelerates the ship either in shallow regions or when a mammal is close to the surface. In the deep-water region, the noise-dominant profile remains consistently at 16 kt, except at a distance of approximately 275 NM, where a mammal is present near the surface. This can be attributed to the low noise levels received by mammals in deep waters due to the upward refracting sound-speed profile. Hence, slowing down is deemed unnecessary in such circumstances. When comparing Case B1 and B2, the noise-dominant profile is nearly identical, with the difference observed when the speed decreases to 8 kt at approximately 450 NM in the case of B2, where bathymetry is taken into account. Fig. \ref{fig_Case_B} illustrates that at this location, the bathymetry undergoes a transition from deep to shallow, having an ascending structure that redirects noise signals toward the mammal. This could be a plausible explanation for this observation. These results further enhance MOOF's adaptive capabilities, demonstrating that the framework prevents unnecessary slowdowns depending on the environment.

A quantitative comparison of Case B1 and B2 in Table \ref{tab_results} shows minimal differences. Case B2, a realistic deep-water case, demonstrates a reduction of 0.95 dB in total noise intensity (88.9\% lower compared to the fuel-dominant solution) 
through the optimization of sailing speed, requiring an additional 91.52 MT of fuel (0.8\% higher compared to the fuel-dominant solution).

\section{Conclusion} \label{sec_conclusion}

In this paper, we introduced a multi-objective optimization framework to determine the optimal operational ship speed, minimizing underwater radiated noise and total fuel consumption along a fixed ship route, while subjected to voyage constraints. This framework incorporated a 2D underwater environment that considered marine mammals from different audiogram groups randomly dispersed amidst variable bathymetry and a range-independent sound speed profile. An empirical model is used to estimate the near-field source noise level, while a ray-based model serves as the underwater propagation model to calculate received noise levels. Using the received levels, an objective function is defined that encapsulates the acoustic footprints of the voyage. In parallel, another objective function is defined to represent the total fuel consumption, calculated using a ship performance model. The non-dominated sorting genetic algorithm-II is incorporated to determine solutions on the normalized Pareto front. Subsequently, a multi-criteria decision-making method, utilizing TOPSIS, is used to identify the trade-off solution among the normalized Pareto front solutions.

In this study, three main cases were examined: a test case (Case T) having simplified environment configuration to illustrate the environmental influence on the speed profile, and two practical cases representing shallow water (Case A) and deep water (Case B) regions to demonstrate the efficacy of the framework. A  6900 TEU containership sailing in the Salish Sea and the Celtic Sea is simulated corresponding to Case A and Case B, respectively.  A pragmatic route is chosen and discretized into multiple sailing legs. At each leg, the aforementioned acoustic solver is solved for frequencies ranging from $10$ to $10^4$ Hz, as well as the fuel consumption is estimated using the empirical model. Subsequently, NSGA-II is employed to determine three optimal speed profiles - the noise-dominant, fuel-dominant, and trade-off profiles - which are analyzed and compared.

Case T showcased that within the deep-water environment, shadow and convergence zones are formed near the surface. The proposed framework of the ship's speed profile strategically leveraged this phenomenon, allowing for higher speeds within the shadow zone and conversely, adjusting the speed profile elsewhere.
Moreover, MOOF also leveraged the sound refraction phenomenon, causing a peak in the speed profile precisely positioned above the marine mammal. These case studies highlighted that this framework delivers an optimal speed profile adaptive to the environment.

Based on the results of the practical case studies, the trade-off solution of Case A showed a significant reduction of 1.31 dB in the total intensity of URN (95\% lower compared to the fuel-dominant solution), with an increase of 4.33 MT in fuel consumption (0.6\% higher than the fuel-dominant solution). In Case B, similar results were observed, with the total intensity of URN reduced by 0.95 dB (88.9\% lower compared to the fuel-dominant solution) 
and an increase of 91.52 MT of fuel (0.8\% higher compared to the fuel-dominant solution). Hence, quantitatively this framework enables the reduction of total URN intensity levels without significantly compromising fuel consumption in a realistic shipping route, considering some simple assumptions. The optimized trade-off speed profile in all the cases is closely aligned with the fuel-dominant profile while incorporating some characteristics of the noise-dominant profile. Specifically, the speed profile emphasized slowdowns in proximity to Phocid Carnivores due to their higher sensitivity to lower frequencies. This proves the MOOF’s adaptability to various marine mammal types in the environment, prioritizing slowdowns for susceptible mammal types. Furthermore, when comparing Case B to Case A, the framework exhibited versatility in different environments. Specifically, Case B illustrated the framework's ability to prevent unnecessary slowdown, in situations where the received noise level is nearly negligible due to the upward refraction of sound.

Although the proposed framework demonstrates efficacy in mitigating noise pollution, certain limitations potentially affect its practicality. Specifically, the framework assumes a fixed location for marine mammals, without accounting for their mobility or the dynamic variation of environmental conditions, such as weather variations and oceanic currents. Future research could address these limitations by incorporating spatiotemporal models, such as a state-space model for mobility \cite{jonsen2005robust} and a population distribution model \cite{melo2020ecological}. Moreover, ship fuel consumption used in this study can be replaced with complex models that consider factors such as weather, policies, operational conditions, and parameters \cite{fan2022review}. Furthermore, replacing the current acoustic models with data-driven deep-learning models, offering improved accuracy and efficiency. For estimating near-field URN levels, there is potential to replace the near-field empirical model with a statistical regression model from the Enhancing Cetacean Habitat and Observation program \cite{MacGillivray2022}. Moreover, the propagation model can be replaced with a generalized deep-learning model for far-field acoustic transmission loss predictions \cite{deo2024continual, deo2022predicting, Mallik2022, mallik2024deep}. The ray-based model used in this study neglects complex oceanographic conditions and the range-dependent variations of these conditions, which are crucial, especially for real-time ship voyage applications. Machine learning models offer a promising solution to address these challenges, providing faster and more accurate computations for far-field propagation \cite{Brissaud2022}. Furthermore, the MOOF can be extended to a 3D ocean environment with moving marine mammals where dynamic path planning capabilities can be incorporated, alongside the optimization of sailing speed.


\section*{Acknowledgments}
The present study is supported by Mitacs, Transport Canada, and Clear Seas through the Quiet-Vessel Initiative (QVI) program. The authors would like to express their gratitude to Dr. Paul Blomerous and Ms. Tessa Coulthard for their valuable feedback and suggestions. Additionally, we gratefully acknowledge the contributions of Dr. David Rosen and Dr. Andrew W Trites, who provided insights on the marine biology aspect of the study. This research was also supported in part through computational resources and services provided by Advanced Research Computing (ARC) at the University of British Columbia and Compute Canada.

\appendix

\section{Model verification}
\label{sec_lap_verif}
\begin{figure}[ht]
\centering
\includegraphics[width=6.5cm]{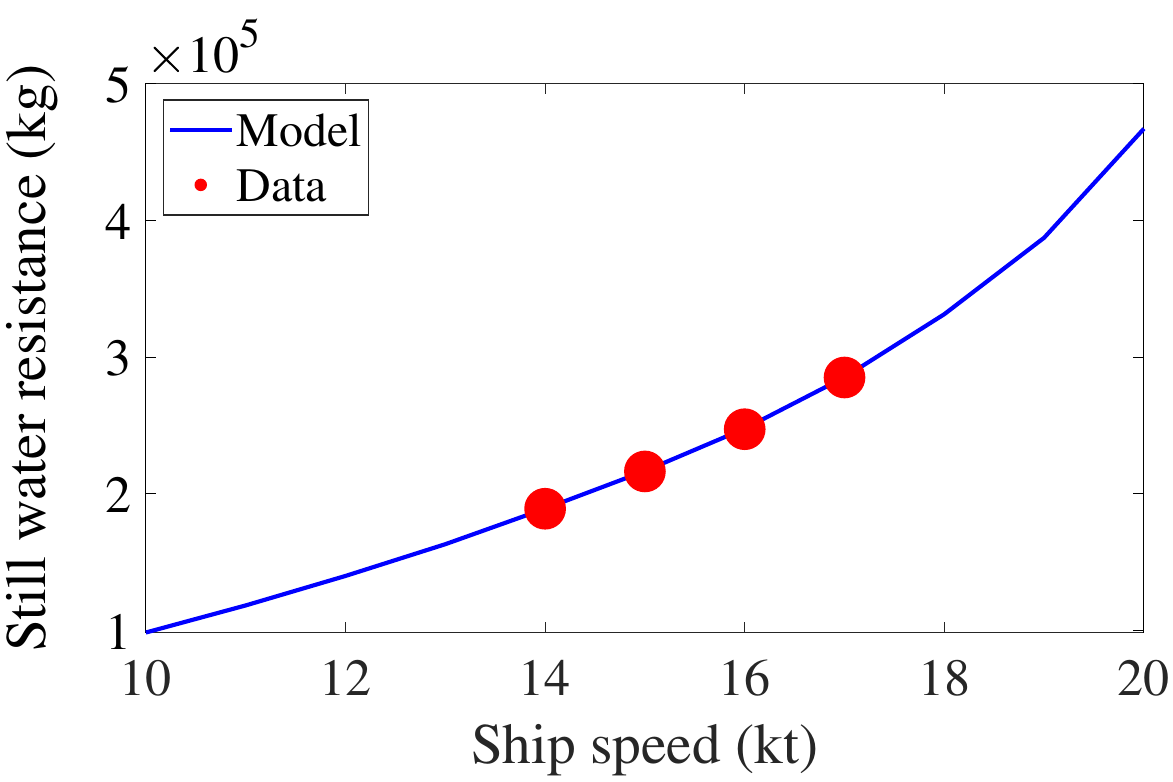}
\caption{Comparison of total water resistance between the model and data from \cite{Keller1973} at different speeds.}
\label{fig_val}
\end{figure}

As mentioned previously, we use the Lap-Keller method to calculate the still water resistance, which enables us to determine the hourly fuel consumption. In this section, the still water resistance is validated with a test case from \cite{Keller1973}. The comparison for various ship speeds is illustrated in Fig. \ref{fig_val}. Detailed information regarding the ship characteristics of the test case can be found in \cite{Keller1973}. It is evident that the model calculations are consistent with the test case data points. The maximum percent error is 0.18\% for still water resistance.

\section{Properties of ship and environment}
\label{sec_props_ship_env}
\begin{figure}[ht]
\centering
\includegraphics[width=6.5cm]{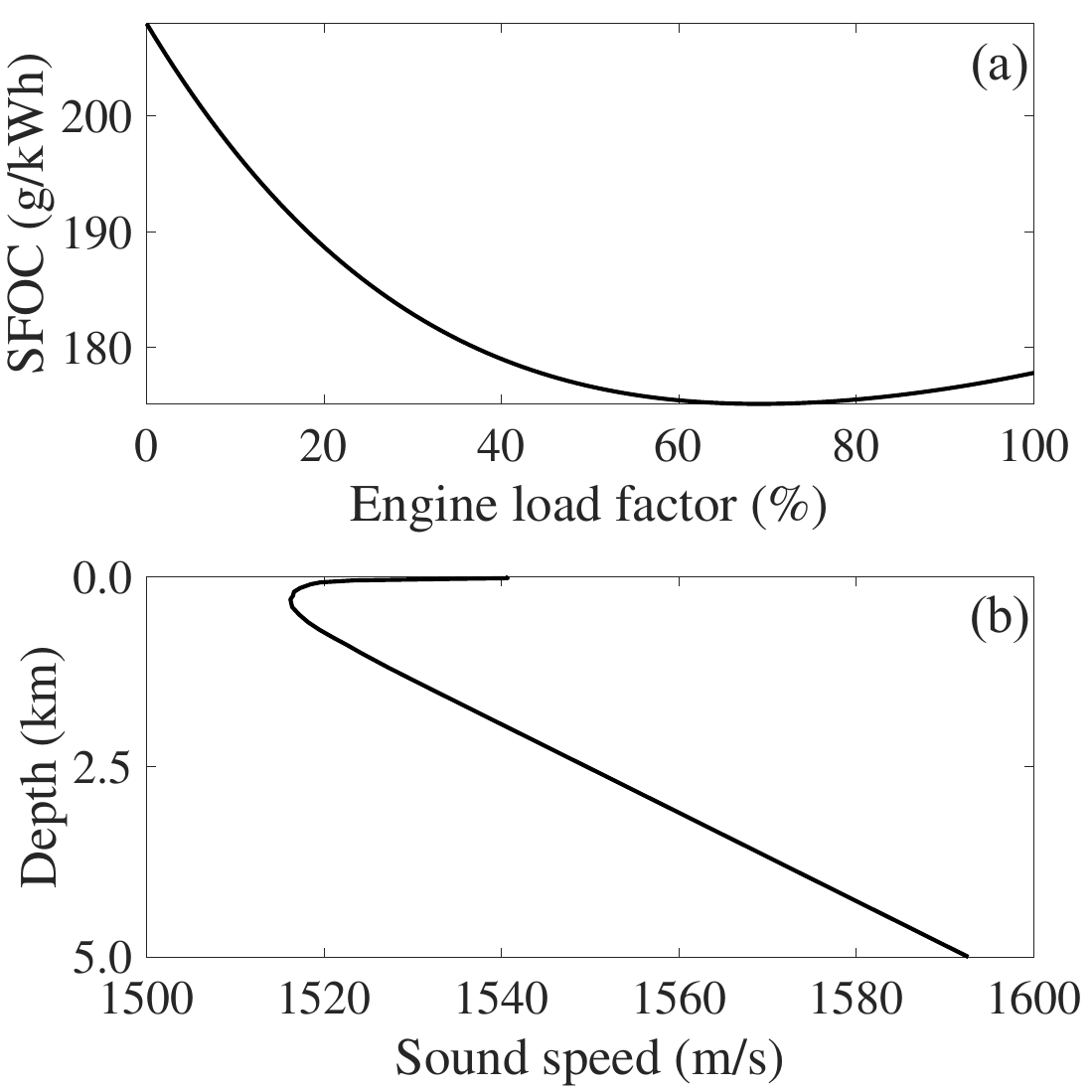}
\caption{(a) SFOC vs engine load percentage \cite{Verma2021}, and (b) SSP \cite{Audet1974AESDSP}.}
\label{fig_sfoc}
\end{figure}

This section delineates the essential properties pertinent to the ship and its operational environment used in the case studies. Initially, for the ship, Table \ref{tab_ship_prop} presents a comprehensive overview of the ship characteristics utilized in the calculations for water resistance. Subsequently, upon determining the resistance, the required propulsion power is computed, along with its associated fuel consumption which requires SFOC. Fig. \ref{fig_sfoc} depicts the SFOC curve of the engine \cite{Tzortzis2021} used in this study.

The environmental conditions comprise factors such as the ship path, bathymetry, sound speed profile, and mammal locations. Fig. \ref{fig_Case_A} and Fig. \ref{fig_Case_B} depict the trajectory of the ship and the position of the marine mammal along with the audiogram group it belongs to in Cases A and B, respectively. The range-independent SSP \cite{Audet1974AESDSP} used in Bellhop is plotted in Fig. \ref{fig_sfoc}.

\begin{table}[ht]
\begin{tabular}{ll}
\hline
Parameters                    & Value      \\ \hline
Length between perpendiculars & 258.4 m    \\
Breadth on waterline          & 42.8 m     \\
Draft                         & 14.55 m    \\
Block coefficient             & 0.6234     \\
Midship coefficient           & 0.9750      \\
Deadweight                    & 72294 tons \\ \hline
\end{tabular}
\caption{Ship characteristics.}
\label{tab_ship_prop}
\end{table}

\section{Parameter tuning of NSGA-II}
\label{sec_hyperparam}

This section presents the parameter tuning of NSGA-II using SMAC3, where the parameters are optimized based on a metric to enhance performance. Before proceeding, consider a general parameterized algorithm, \(\Lambda\), for solving problem instance(s), \(I\). The parameter or configuration space, \(\Theta\), constitutes all possible parameters of \(\Lambda\). Parameter tuning involves finding the optimal parameter \(\theta^* \in \Theta\) that maximizes a performance metric \(m(\theta)\). A general parameter tuning problem is defined as follows:

\begin{equation}
\begin{aligned}
\theta^* &= \max_{\theta \in \Theta} m(\theta), \\
\text{where} \quad m(\theta) &= f(\theta \mid I, P_I, P_{\zeta}, t).
\end{aligned}
\end{equation}

 Here, \(P_I\) represents the distribution over admissible instances of the problem. The function \( \zeta(\theta, i, t) = \zeta(\Lambda(\theta), i, t) \) assigns a cost value to each configuration \( \theta \) when the algorithm \( \Lambda(\theta) \) is executed on instance \( i \) from set \( I \) for a specified runtime \( t \). In a model-based optimization algorithm, this cost is modeled as another distribution function, \( \zeta \sim P_{\zeta}(\zeta \mid \theta, i, t) \).

In this study, NSGA-II is the parameterized algorithm where the parameter space is defined by population size, crossover probability, and mutation probability. The bounds for this space are set according to Table \ref{tab_config}. The number of generations for each configuration is capped at 500 to avoid excessive computational expense. The random forest model in SMAC3 is used as the surrogate model.

\begin{table}[]
\begin{tabular}{llll}
\hline
Parameter                                                       & Type       & Range          & \begin{tabular}[c]{@{}l@{}}Optimal\\ value\end{tabular} \\ \hline
\begin{tabular}[c]{@{}l@{}}Population\\ size\end{tabular}       & Integer    & {[}50,250{]}   & 200                                                     \\
\begin{tabular}[c]{@{}l@{}}Crossover\\ probability\end{tabular} & Continuous & {[}0.01,1.0{]} & 0.88                                                     \\
\begin{tabular}[c]{@{}l@{}}Mutation\\ probability\end{tabular}  & Continuous & {[}0.01,1.0{]} & 0.025                                                     \\ \hline
\end{tabular}
\caption{The parameter space and optimal parameters of NSGA-II.}
\label{tab_config}
\end{table}

The hypervolume indicator \cite{shang2020survey} is commonly utilized to assess the quality of Pareto solution sets from evolutionary multi-objective optimization. Given a point set \( A \subset \mathbb{R}^m \) and a reference point \( r \in \mathbb{R}^m \), the hypervolume of the point set \( A \) is defined as:

\begin{equation}
    \text{HV}(A, r) = \mathcal{L} \left( \bigcup_{a \in A} \{ b \mid a \geq b \geq r \} \right),
\end{equation}
where \( \mathcal{L} \) denotes the Lebesgue measure of a set on $\mathbb{R}^m$. The notation \( a \geq b \) means that \( a \) dominates \( b \) (i.e., \( a_i \geq b_i \) for all \( i = 1, \ldots, m \) and \( a_{j} > b_{j} \) for at least one \( j = 1, \ldots, m \) in the maximization case) \cite{shang2020survey}. This quantity is the performance metric, \(m(\theta)\), used in this study, which is maximized during the tuning process.

\begin{figure}[ht]
\centering
\includegraphics[width=7cm]{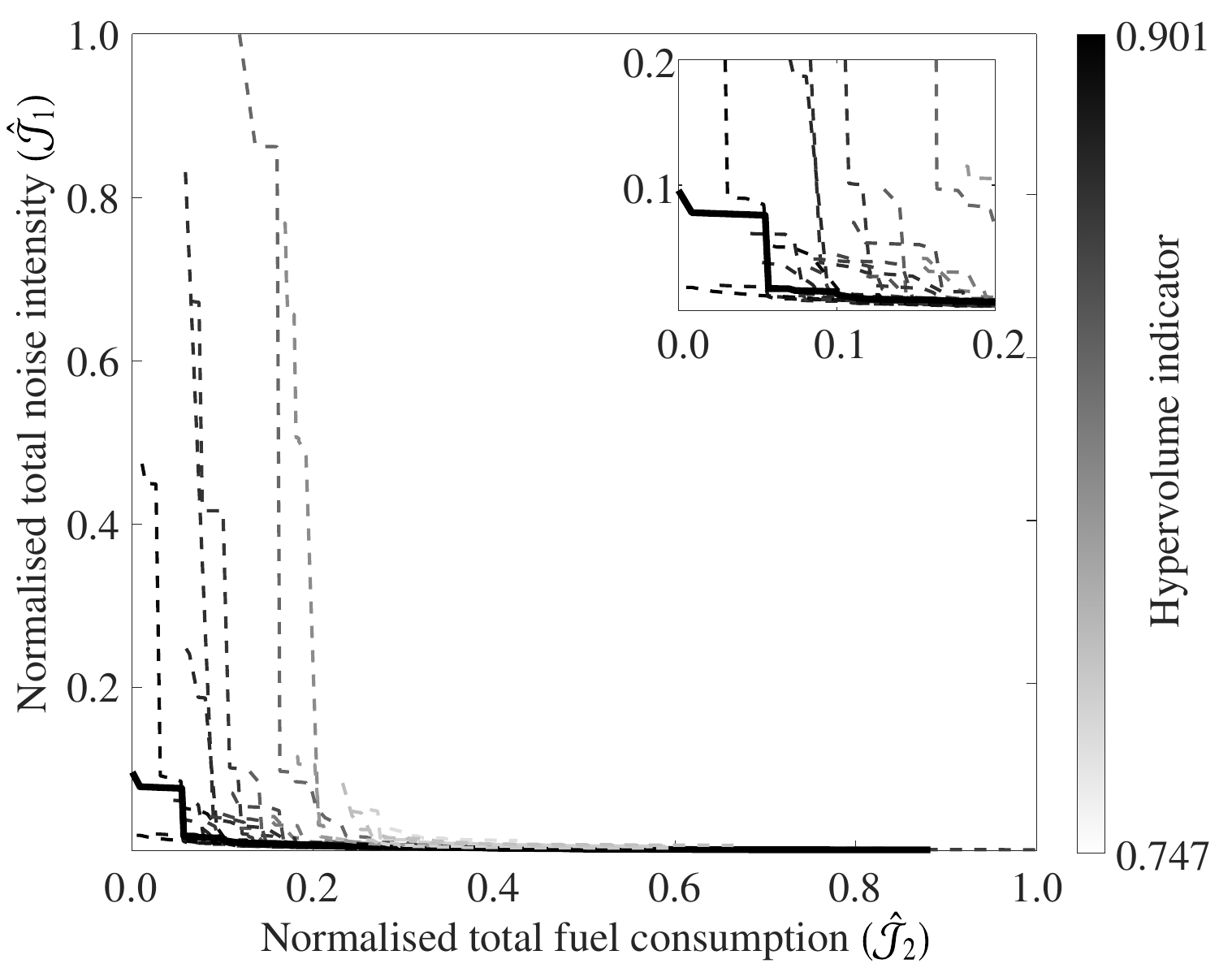}
\caption{Pareto front of the last 25 trials with varying configurations from SMAC3, color-coded by the hypervolume indicator. The optimal Pareto front is represented by the solid line, with its configuration detailed in Table \ref{tab_config}.}
\label{fig_hyp}
\end{figure}

SMAC3 is used to find the optimal parameters of NSGA-II by maximizing the hypervolume indicator for a particular problem instance, i.e., Case A2. This case is specifically chosen for tuning as it most closely resembles a realistic case study. Fig. \ref{fig_hyp} shows the Pareto front of these trials, with the Pareto front having the maximum hypervolume indicated by a solid line and the corresponding optimal parameters are tabulated in Table \ref{tab_config}.

 \bibliographystyle{elsarticle-num-names} 
 \bibliography{biblograph}

\begin{thebibliography}{70}
\expandafter\ifx\csname natexlab\endcsname\relax\def\natexlab#1{#1}\fi
\providecommand{\url}[1]{\texttt{#1}}
\providecommand{\href}[2]{#2}
\providecommand{\path}[1]{#1}
\providecommand{\DOIprefix}{doi:}
\providecommand{\ArXivprefix}{arXiv:}
\providecommand{\URLprefix}{URL: }
\providecommand{\Pubmedprefix}{pmid:}
\providecommand{\doi}[1]{\href{http://dx.doi.org/#1}{\path{#1}}}
\providecommand{\Pubmed}[1]{\href{pmid:#1}{\path{#1}}}
\providecommand{\bibinfo}[2]{#2}
\ifx\xfnm\relax \def\xfnm[#1]{\unskip,\space#1}\fi
\bibitem[{Kaplan and Solomon(2016)}]{Kaplan2016}
\bibinfo{author}{M.~B. Kaplan}, \bibinfo{author}{S.~Solomon},
\newblock \bibinfo{title}{{A coming boom in commercial shipping? The potential for rapid growth of noise from commercial ships by 2030}},
\newblock \bibinfo{journal}{Marine Policy} \bibinfo{volume}{73} (\bibinfo{year}{2016}) \bibinfo{pages}{119--121}.
\bibitem[{Jägerbrand et~al.(2019)Jägerbrand, Brutemark, {Barthel Svedén}, and Gren}]{Jagerbrand2019}
\bibinfo{author}{A.~K. Jägerbrand}, \bibinfo{author}{A.~Brutemark}, \bibinfo{author}{J.~{Barthel Svedén}}, \bibinfo{author}{I.-M. Gren},
\newblock \bibinfo{title}{A review on the environmental impacts of shipping on aquatic and nearshore ecosystems},
\newblock \bibinfo{journal}{Science of The Total Environment} \bibinfo{volume}{695} (\bibinfo{year}{2019}) \bibinfo{pages}{133637}.
\bibitem[{Peng et~al.(2015)Peng, Zhao, and Liu}]{Peng2015}
\bibinfo{author}{C.~Peng}, \bibinfo{author}{X.~Zhao}, \bibinfo{author}{G.~Liu},
\newblock \bibinfo{title}{{Noise in the sea and its impacts on marine organisms}},
\newblock \bibinfo{journal}{International Journal of Environmental Research and Public Health} \bibinfo{volume}{12} (\bibinfo{year}{2015}) \bibinfo{pages}{12304--12323}.
\bibitem[{Miksis-Olds and Nichols(2016)}]{MiksisOlds2016}
\bibinfo{author}{J.~L. Miksis-Olds}, \bibinfo{author}{S.~M. Nichols},
\newblock \bibinfo{title}{{Is low frequency ocean sound increasing globally?}},
\newblock \bibinfo{journal}{The Journal of the Acoustical Society of America} \bibinfo{volume}{139} (\bibinfo{year}{2016}) \bibinfo{pages}{501--511}.
\bibitem[{McDonald et~al.(2006)McDonald, Hildebrand, and Wiggins}]{McDonald2006}
\bibinfo{author}{M.~A. McDonald}, \bibinfo{author}{J.~A. Hildebrand}, \bibinfo{author}{S.~M. Wiggins},
\newblock \bibinfo{title}{{Increases in deep ocean ambient noise in the Northeast Pacific west of San Nicolas Island, California}},
\newblock \bibinfo{journal}{The Journal of the Acoustical Society of America} \bibinfo{volume}{120} (\bibinfo{year}{2006}) \bibinfo{pages}{711--718}.
\bibitem[{Erbe et~al.(2019)Erbe, Marley, Schoeman, Smith, Trigg, and Embling}]{Erbe2019}
\bibinfo{author}{C.~Erbe}, \bibinfo{author}{S.~A. Marley}, \bibinfo{author}{R.~P. Schoeman}, \bibinfo{author}{J.~N. Smith}, \bibinfo{author}{L.~E. Trigg}, \bibinfo{author}{C.~B. Embling},
\newblock \bibinfo{title}{{The Effects of Ship Noise on Marine Mammals—A Review}},
\newblock \bibinfo{journal}{Frontiers in Marine Science} \bibinfo{volume}{6} (\bibinfo{year}{2019}).
\bibitem[{Richardson et~al.(1995)Richardson, Greene, Malme, and Thomson}]{Richardson1995}
\bibinfo{author}{W.~J. Richardson}, \bibinfo{author}{C.~R. Greene}, \bibinfo{author}{C.~I. Malme}, \bibinfo{author}{D.~H. Thomson},
\newblock \bibinfo{title}{{Significance of Responses and Noise Impacts}},
\newblock \bibinfo{journal}{Marine Mammals and Noise}  (\bibinfo{year}{1995}) \bibinfo{pages}{387--424}.
\bibitem[{Gomez et~al.(2016)Gomez, Lawson, Wright, Buren, Tollit, and Lesage}]{Gomez2016}
\bibinfo{author}{C.~Gomez}, \bibinfo{author}{J.~Lawson}, \bibinfo{author}{A.~Wright}, \bibinfo{author}{A.~Buren}, \bibinfo{author}{D.~Tollit}, \bibinfo{author}{V.~Lesage},
\newblock \bibinfo{title}{A systematic review on the behavioural responses of wild marine mammals to noise: the disparity between science and policy},
\newblock \bibinfo{journal}{Canadian Journal of Zoology} \bibinfo{volume}{94} (\bibinfo{year}{2016}) \bibinfo{pages}{801--819}.
\bibitem[{Wright et~al.(2007)Wright, Soto, Baldwin, Bateson, Beale, Clark, Deak, Edwards, Fern{\'a}ndez, Godinho et~al.}]{Wright2007}
\bibinfo{author}{A.~J. Wright}, \bibinfo{author}{N.~A. Soto}, \bibinfo{author}{A.~L. Baldwin}, \bibinfo{author}{M.~Bateson}, \bibinfo{author}{C.~M. Beale}, \bibinfo{author}{C.~Clark}, \bibinfo{author}{T.~Deak}, \bibinfo{author}{E.~F. Edwards}, \bibinfo{author}{A.~Fern{\'a}ndez}, \bibinfo{author}{A.~Godinho}, et~al.,
\newblock \bibinfo{title}{Do marine mammals experience stress related to anthropogenic noise?},
\newblock \bibinfo{journal}{International Journal of Comparative Psychology} \bibinfo{volume}{20} (\bibinfo{year}{2007}).
\bibitem[{Southall et~al.(2019)Southall, Finneran, Reichmuth, Nachtigall, Ketten, Bowles, Ellison, Nowacek, and Tyack}]{Southall2019}
\bibinfo{author}{E.~B.~L. Southall}, \bibinfo{author}{J.~J. Finneran}, \bibinfo{author}{C.~Reichmuth}, \bibinfo{author}{P.~E. Nachtigall}, \bibinfo{author}{D.~R. Ketten}, \bibinfo{author}{A.~E. Bowles}, \bibinfo{author}{W.~T. Ellison}, \bibinfo{author}{D.~P. Nowacek}, \bibinfo{author}{P.~L. Tyack},
\newblock \bibinfo{title}{{Marine mammal noise exposure criteria: Updated scientific recommendations for residual hearing effects}},
\newblock \bibinfo{journal}{Aquatic Mammals} \bibinfo{volume}{45} (\bibinfo{year}{2019}) \bibinfo{pages}{125--232}.
\bibitem[{Chou et~al.(2021)Chou, Southall, Robards, and Rosenbaum}]{Chou2021}
\bibinfo{author}{E.~Chou}, \bibinfo{author}{B.~L. Southall}, \bibinfo{author}{M.~D. Robards}, \bibinfo{author}{H.~C. Rosenbaum},
\newblock \bibinfo{title}{International policy, recommendations, actions and mitigation efforts of anthropogenic underwater noise},
\newblock \bibinfo{journal}{Ocean \& Coastal Management} \bibinfo{volume}{202} (\bibinfo{year}{2021}) \bibinfo{pages}{105427--}.
\bibitem[{Leaper(2019)}]{Leaper2019}
\bibinfo{author}{R.~Leaper},
\newblock \bibinfo{title}{The role of slower vessel speeds in reducing greenhouse gas emissions, underwater noise and collision risk to whales},
\newblock \bibinfo{journal}{Frontiers in Marine Science} \bibinfo{volume}{6} (\bibinfo{year}{2019}).
\bibitem[{Joy et~al.(2019)Joy, Tollit, Wood, MacGillivray, Li, Trounce, and Robinson}]{Joy2019}
\bibinfo{author}{R.~Joy}, \bibinfo{author}{D.~Tollit}, \bibinfo{author}{J.~Wood}, \bibinfo{author}{A.~MacGillivray}, \bibinfo{author}{Z.~Li}, \bibinfo{author}{K.~Trounce}, \bibinfo{author}{O.~Robinson},
\newblock \bibinfo{title}{Potential benefits of vessel slowdowns on endangered southern resident killer whales},
\newblock \bibinfo{journal}{Frontiers in Marine Science} \bibinfo{volume}{6} (\bibinfo{year}{2019}).
\bibitem[{Authority(2022{\natexlab{a}})}]{ECHO_boundarypass}
\bibinfo{author}{V.~F.~P. Authority},
\newblock \bibinfo{title}{2023 haro strait and boundary pass voluntary ship slowdown},
\newblock \bibinfo{journal}{ECHO Program}  (\bibinfo{year}{2022}{\natexlab{a}}).
\bibitem[{Authority(2022{\natexlab{b}})}]{ECHO_swiftsure}
\bibinfo{author}{V.~F.~P. Authority},
\newblock \bibinfo{title}{2023 swiftsure bank voluntary ship slowdown},
\newblock \bibinfo{journal}{ECHO Program}  (\bibinfo{year}{2022}{\natexlab{b}}).
\bibitem[{Psaraftis and Kontovas(2014)}]{Psaraftis2014}
\bibinfo{author}{H.~N. Psaraftis}, \bibinfo{author}{C.~A. Kontovas},
\newblock \bibinfo{title}{Ship speed optimization: Concepts, models and combined speed-routing scenarios},
\newblock \bibinfo{journal}{Transportation Research Part C: Emerging Technologies} \bibinfo{volume}{44} (\bibinfo{year}{2014}) \bibinfo{pages}{52--69}.
\bibitem[{Psaraftis and Kontovas(2013)}]{Psaraftis2013}
\bibinfo{author}{H.~N. Psaraftis}, \bibinfo{author}{C.~A. Kontovas},
\newblock \bibinfo{title}{Speed models for energy-efficient maritime transportation: A taxonomy and survey},
\newblock \bibinfo{journal}{Transportation Research Part C: Emerging Technologies} \bibinfo{volume}{26} (\bibinfo{year}{2013}) \bibinfo{pages}{331--351}.
\bibitem[{Li et~al.(2020)Li, Sun, Guo, Du, and Li}]{Li2020}
\bibinfo{author}{X.~Li}, \bibinfo{author}{B.~Sun}, \bibinfo{author}{C.~Guo}, \bibinfo{author}{W.~Du}, \bibinfo{author}{Y.~Li},
\newblock \bibinfo{title}{Speed optimization of a container ship on a given route considering voluntary speed loss and emissions},
\newblock \bibinfo{journal}{Applied Ocean Research} \bibinfo{volume}{94} (\bibinfo{year}{2020}) \bibinfo{pages}{101995}.
\bibitem[{Tzortzis and Sakalis(2021)}]{Tzortzis2021}
\bibinfo{author}{G.~Tzortzis}, \bibinfo{author}{G.~Sakalis},
\newblock \bibinfo{title}{A dynamic ship speed optimization method with time horizon segmentation},
\newblock \bibinfo{journal}{Ocean Engineering} \bibinfo{volume}{226} (\bibinfo{year}{2021}) \bibinfo{pages}{108840}.
\bibitem[{Yu et~al.(2021)Yu, Fang, Fu, Liu, and Chen}]{Yu2021}
\bibinfo{author}{H.~Yu}, \bibinfo{author}{Z.~Fang}, \bibinfo{author}{X.~Fu}, \bibinfo{author}{J.~Liu}, \bibinfo{author}{J.~Chen},
\newblock \bibinfo{title}{Literature review on emission control-based ship voyage optimization},
\newblock \bibinfo{journal}{Transportation Research Part D: Transport and Environment} \bibinfo{volume}{93} (\bibinfo{year}{2021}) \bibinfo{pages}{102768}.
\bibitem[{Ma et~al.(2020)Ma, Ma, Jin, and Ma}]{Ma2020}
\bibinfo{author}{D.~Ma}, \bibinfo{author}{W.~Ma}, \bibinfo{author}{S.~Jin}, \bibinfo{author}{X.~Ma},
\newblock \bibinfo{title}{Method for simultaneously optimizing ship route and speed with emission control areas},
\newblock \bibinfo{journal}{Ocean Engineering} \bibinfo{volume}{202} (\bibinfo{year}{2020}) \bibinfo{pages}{107170}.
\bibitem[{Ma et~al.(2021)Ma, Ma, Ma, Zhang, and Wang}]{Ma2021}
\bibinfo{author}{W.~Ma}, \bibinfo{author}{D.~Ma}, \bibinfo{author}{Y.~Ma}, \bibinfo{author}{J.~Zhang}, \bibinfo{author}{D.~Wang},
\newblock \bibinfo{title}{Green maritime: a routing and speed multi-objective optimization strategy},
\newblock \bibinfo{journal}{Journal of Cleaner Production} \bibinfo{volume}{305} (\bibinfo{year}{2021}) \bibinfo{pages}{127179}.
\bibitem[{Wang et~al.(2019)Wang, Mao, and Eriksson}]{Wang2019}
\bibinfo{author}{H.~Wang}, \bibinfo{author}{W.~Mao}, \bibinfo{author}{L.~Eriksson},
\newblock \bibinfo{title}{A three-dimensional dijkstra's algorithm for multi-objective ship voyage optimization},
\newblock \bibinfo{journal}{Ocean Engineering} \bibinfo{volume}{186} (\bibinfo{year}{2019}) \bibinfo{pages}{106131}.
\bibitem[{Wen et~al.(2017)Wen, Pacino, Kontovas, and Psaraftis}]{Wen2017}
\bibinfo{author}{M.~Wen}, \bibinfo{author}{D.~Pacino}, \bibinfo{author}{C.~Kontovas}, \bibinfo{author}{H.~Psaraftis},
\newblock \bibinfo{title}{A multiple ship routing and speed optimization problem under time, cost and environmental objectives},
\newblock \bibinfo{journal}{Transportation Research Part D: Transport and Environment} \bibinfo{volume}{52} (\bibinfo{year}{2017}) \bibinfo{pages}{303--321}.
\bibitem[{Khatami et~al.(2023)Khatami, Chen, and Chen}]{Khatami2023}
\bibinfo{author}{R.~Khatami}, \bibinfo{author}{B.~Chen}, \bibinfo{author}{Y.~C. Chen},
\newblock \bibinfo{title}{Optimal voyage scheduling of all-electric ships considering underwater radiated noise},
\newblock \bibinfo{journal}{Transportation Research Part C: Emerging Technologies} \bibinfo{volume}{148} (\bibinfo{year}{2023}) \bibinfo{pages}{104024}.
\bibitem[{Stojanovic and Preisig(2009)}]{Stojanovic2009}
\bibinfo{author}{M.~Stojanovic}, \bibinfo{author}{J.~Preisig},
\newblock \bibinfo{title}{Underwater acoustic communication channels: Propagation models and statistical characterization},
\newblock \bibinfo{journal}{IEEE Communications Magazine} \bibinfo{volume}{47} (\bibinfo{year}{2009}) \bibinfo{pages}{84--89}.
\bibitem[{DeSanto(1979)}]{desanto1979theoretical}
\bibinfo{author}{J.~A. DeSanto},
\newblock \bibinfo{title}{Theoretical methods in ocean acoustics},
\newblock in: \bibinfo{booktitle}{Ocean Acoustics}, \bibinfo{publisher}{Springer}, \bibinfo{year}{1979}, pp. \bibinfo{pages}{7--77}.
\bibitem[{Oliveira et~al.(2021)Oliveira, Lin, and Porter}]{oliveira2021underwater}
\bibinfo{author}{T.~C. Oliveira}, \bibinfo{author}{Y.-T. Lin}, \bibinfo{author}{M.~B. Porter},
\newblock \bibinfo{title}{Underwater sound propagation modeling in a complex shallow water environment},
\newblock \bibinfo{journal}{Frontiers in Marine Science} \bibinfo{volume}{8} (\bibinfo{year}{2021}) \bibinfo{pages}{751327}.
\bibitem[{Jensen et~al.(2011)Jensen, Kuperman, Porter, Schmidt, and Tolstoy}]{jensen2011computational}
\bibinfo{author}{F.~B. Jensen}, \bibinfo{author}{W.~A. Kuperman}, \bibinfo{author}{M.~B. Porter}, \bibinfo{author}{H.~Schmidt}, \bibinfo{author}{A.~Tolstoy}, \bibinfo{title}{Computational ocean acoustics}, volume \bibinfo{volume}{2011}, \bibinfo{publisher}{Springer}, \bibinfo{year}{2011}.
\bibitem[{Etter(2018)}]{etter2018underwater}
\bibinfo{author}{P.~C. Etter}, \bibinfo{title}{Underwater acoustic modeling and simulation}, \bibinfo{publisher}{CRC press}, \bibinfo{year}{2018}.
\bibitem[{Porter and Bucker(1987)}]{porter1987gaussian}
\bibinfo{author}{M.~B. Porter}, \bibinfo{author}{H.~P. Bucker},
\newblock \bibinfo{title}{Gaussian beam tracing for computing ocean acoustic fields},
\newblock \bibinfo{journal}{The Journal of the Acoustical Society of America} \bibinfo{volume}{82} (\bibinfo{year}{1987}) \bibinfo{pages}{1349--1359}.
\bibitem[{Porter(2011)}]{porter2011bellhop}
\bibinfo{author}{M.~B. Porter},
\newblock \bibinfo{title}{The bellhop manual and user’s guide: Preliminary draft},
\newblock \bibinfo{journal}{Heat, Light, and Sound Research, Inc., La Jolla, CA, USA, Tech. Rep} \bibinfo{volume}{260} (\bibinfo{year}{2011}).
\bibitem[{Silva et~al.(2016)Silva, Neves, and Horta}]{Silva2016}
\bibinfo{author}{A.~D. Silva}, \bibinfo{author}{R.~F. Neves}, \bibinfo{author}{N.~Horta}, \bibinfo{title}{{Multi-objective optimization}}, \bibinfo{year}{2016}.
\bibitem[{Kuperman(2003)}]{KUPERMAN2003317}
\bibinfo{author}{W.~A. Kuperman},
\newblock \bibinfo{title}{Underwater acoustics},
\newblock in: \bibinfo{editor}{R.~A. Meyers} (Ed.), \bibinfo{booktitle}{Encyclopedia of Physical Science and Technology (Third Edition)}, \bibinfo{edition}{third edition} ed., \bibinfo{publisher}{Academic Press}, \bibinfo{address}{New York}, \bibinfo{year}{2003}, pp. \bibinfo{pages}{317--338}.
\bibitem[{Zhu et~al.(2022)Zhu, Gaggero, Makris, and Ratilal}]{Zhu2022}
\bibinfo{author}{C.~Zhu}, \bibinfo{author}{T.~Gaggero}, \bibinfo{author}{N.~C. Makris}, \bibinfo{author}{P.~Ratilal},
\newblock \bibinfo{title}{{Underwater Sound Characteristics of a Ship with Controllable Pitch Propeller}},
\newblock \bibinfo{journal}{Journal of Marine Science and Engineering} \bibinfo{volume}{10} (\bibinfo{year}{2022}).
\bibitem[{Wittekind(2014)}]{Wittekind2014}
\bibinfo{author}{D.~Wittekind},
\newblock \bibinfo{title}{A simple model for the underwater noise source level of ships},
\newblock \bibinfo{journal}{Journal of ship production and design} \bibinfo{volume}{30} (\bibinfo{year}{2014}) \bibinfo{pages}{7--14}.
\bibitem[{Ross(1976)}]{Ross1976}
\bibinfo{author}{D.~Ross}, \bibinfo{title}{{Mechanics of Underwater Noise}}, \bibinfo{year}{1976}.
\bibitem[{Chion et~al.(2019)Chion, Lagrois, and Dupras}]{chion2019meta}
\bibinfo{author}{C.~Chion}, \bibinfo{author}{D.~Lagrois}, \bibinfo{author}{J.~Dupras},
\newblock \bibinfo{title}{A meta-analysis to understand the variability in reported source levels of noise radiated by ships from opportunistic studies},
\newblock \bibinfo{journal}{Frontiers in Marine Science} \bibinfo{volume}{6} (\bibinfo{year}{2019}) \bibinfo{pages}{714}.
\bibitem[{Mallik et~al.(2024)Mallik, Jaiman, and Jelovica}]{mallik2024deep}
\bibinfo{author}{W.~Mallik}, \bibinfo{author}{R.~Jaiman}, \bibinfo{author}{J.~Jelovica},
\newblock \bibinfo{title}{Deep neural network for learning wave scattering and interference of underwater acoustics},
\newblock \bibinfo{journal}{Physics of Fluids} \bibinfo{volume}{36} (\bibinfo{year}{2024}).
\bibitem[{Mallik et~al.(2022)Mallik, Jaiman, and Jelovica}]{Mallik2022}
\bibinfo{author}{W.~Mallik}, \bibinfo{author}{R.~K. Jaiman}, \bibinfo{author}{J.~Jelovica},
\newblock \bibinfo{title}{{Predicting transmission loss in underwater acoustics using convolutional recurrent autoencoder network}},
\newblock \bibinfo{journal}{The Journal of the Acoustical Society of America} \bibinfo{volume}{152} (\bibinfo{year}{2022}) \bibinfo{pages}{1627--1638}.
\bibitem[{Cerveny(1987)}]{cerveny1987ray}
\bibinfo{author}{V.~Cerveny},
\newblock \bibinfo{title}{Ray tracing algorithms in three-dimensional laterally varying layered structures},
\newblock \bibinfo{journal}{Seismic tomography} \bibinfo{volume}{5} (\bibinfo{year}{1987}) \bibinfo{pages}{99--133}.
\bibitem[{Weilgart(2007)}]{Weilgart2007}
\bibinfo{author}{L.~S. Weilgart},
\newblock \bibinfo{title}{{A Brief Review of Known Effects of Noise on Marine Mammals}},
\newblock \bibinfo{journal}{International Journal of Comparative Psychology} \bibinfo{volume}{20} (\bibinfo{year}{2007}).
\bibitem[{Keller(1973)}]{Keller1973}
\bibinfo{author}{W.~H.~m. Keller},
\newblock \bibinfo{title}{{Extended Diagrams for Determining the Resistance and Required Power for Single-Screw Ships.}},
\newblock \bibinfo{journal}{International Shipbuilding Progress} \bibinfo{volume}{20} (\bibinfo{year}{1973}) \bibinfo{pages}{133--142}.
\bibitem[{Lu et~al.(2015)Lu, Turan, Boulougouris, Banks, and Incecik}]{lu2015semi}
\bibinfo{author}{R.~Lu}, \bibinfo{author}{O.~Turan}, \bibinfo{author}{E.~Boulougouris}, \bibinfo{author}{C.~Banks}, \bibinfo{author}{A.~Incecik},
\newblock \bibinfo{title}{A semi-empirical ship operational performance prediction model for voyage optimization towards energy efficient shipping},
\newblock \bibinfo{journal}{Ocean Engineering} \bibinfo{volume}{110} (\bibinfo{year}{2015}) \bibinfo{pages}{18--28}.
\bibitem[{Farag and {\"O}l{\c{c}}er(2020)}]{farag2020development}
\bibinfo{author}{Y.~B. Farag}, \bibinfo{author}{A.~I. {\"O}l{\c{c}}er},
\newblock \bibinfo{title}{The development of a ship performance model in varying operating conditions based on ann and regression techniques},
\newblock \bibinfo{journal}{Ocean Engineering} \bibinfo{volume}{198} (\bibinfo{year}{2020}) \bibinfo{pages}{106972}.
\bibitem[{Yang et~al.(2020)Yang, Chen, Zhao, and Rytter}]{yang2020ship}
\bibinfo{author}{L.~Yang}, \bibinfo{author}{G.~Chen}, \bibinfo{author}{J.~Zhao}, \bibinfo{author}{N.~G.~M. Rytter},
\newblock \bibinfo{title}{Ship speed optimization considering ocean currents to enhance environmental sustainability in maritime shipping},
\newblock \bibinfo{journal}{Sustainability} \bibinfo{volume}{12} (\bibinfo{year}{2020}) \bibinfo{pages}{3649}.
\bibitem[{Wang et~al.(2020)Wang, Lang, Mao, Zhang, and Storhaug}]{wang2020effectiveness}
\bibinfo{author}{H.~Wang}, \bibinfo{author}{X.~Lang}, \bibinfo{author}{W.~Mao}, \bibinfo{author}{D.~Zhang}, \bibinfo{author}{G.~Storhaug},
\newblock \bibinfo{title}{Effectiveness of 2d optimization algorithms considering voluntary speed reduction under uncertain metocean conditions},
\newblock \bibinfo{journal}{Ocean Engineering} \bibinfo{volume}{200} (\bibinfo{year}{2020}) \bibinfo{pages}{107063}.
\bibitem[{van Dooren et~al.(2023)van Dooren, Duhr, and Onder}]{van2023convex}
\bibinfo{author}{S.~van Dooren}, \bibinfo{author}{P.~Duhr}, \bibinfo{author}{C.~H. Onder},
\newblock \bibinfo{title}{Convex modelling for ship speed optimisation},
\newblock \bibinfo{journal}{Ocean Engineering} \bibinfo{volume}{288} (\bibinfo{year}{2023}) \bibinfo{pages}{115947}.
\bibitem[{Deb and Kalyanmoy(2001)}]{Kalyan2001}
\bibinfo{author}{K.~Deb}, \bibinfo{author}{D.~Kalyanmoy}, \bibinfo{title}{Multi-Objective Optimization Using Evolutionary Algorithms}, \bibinfo{publisher}{John Wiley \& Sons, Inc.}, \bibinfo{address}{USA}, \bibinfo{year}{2001}.
\bibitem[{Beheshti and Shamsuddin(2013)}]{Beheshti2013}
\bibinfo{author}{Z.~Beheshti}, \bibinfo{author}{S.~M. Shamsuddin},
\newblock \bibinfo{title}{A review of population-based meta-heuristic algorithm},
\newblock \bibinfo{journal}{International Journal of Advances in Soft Computing and Its Applications} \bibinfo{volume}{5} (\bibinfo{year}{2013}) \bibinfo{pages}{1--35}.
\bibitem[{Verma et~al.(2021)Verma, Pant, and Snasel}]{Verma2021}
\bibinfo{author}{S.~Verma}, \bibinfo{author}{M.~Pant}, \bibinfo{author}{V.~Snasel},
\newblock \bibinfo{title}{{A Comprehensive Review on NSGA-II for Multi-Objective Combinatorial Optimization Problems}},
\newblock \bibinfo{journal}{IEEE Access} \bibinfo{volume}{9} (\bibinfo{year}{2021}) \bibinfo{pages}{57757--57791}.
\bibitem[{Deb and Jain(2013)}]{deb2013evolutionary}
\bibinfo{author}{K.~Deb}, \bibinfo{author}{H.~Jain},
\newblock \bibinfo{title}{An evolutionary many-objective optimization algorithm using reference-point-based nondominated sorting approach, part i: solving problems with box constraints},
\newblock \bibinfo{journal}{IEEE transactions on evolutionary computation} \bibinfo{volume}{18} (\bibinfo{year}{2013}) \bibinfo{pages}{577--601}.
\bibitem[{Goldberg and Deb(1991)}]{goldberg1991comparative}
\bibinfo{author}{D.~E. Goldberg}, \bibinfo{author}{K.~Deb},
\newblock \bibinfo{title}{A comparative analysis of selection schemes used in genetic algorithms},
\newblock in: \bibinfo{booktitle}{Foundations of genetic algorithms}, volume~\bibinfo{volume}{1}, \bibinfo{publisher}{Elsevier}, \bibinfo{year}{1991}, pp. \bibinfo{pages}{69--93}.
\bibitem[{Umbarkar and Sheth(2015)}]{umbarkar2015crossover}
\bibinfo{author}{A.~J. Umbarkar}, \bibinfo{author}{P.~D. Sheth},
\newblock \bibinfo{title}{Crossover operators in genetic algorithms: a review.},
\newblock \bibinfo{journal}{ICTACT journal on soft computing} \bibinfo{volume}{6} (\bibinfo{year}{2015}).
\bibitem[{{Blank} and {Deb}(2020)}]{pymoo}
\bibinfo{author}{J.~{Blank}}, \bibinfo{author}{K.~{Deb}},
\newblock \bibinfo{title}{pymoo: Multi-objective optimization in python},
\newblock \bibinfo{journal}{IEEE Access} \bibinfo{volume}{8} (\bibinfo{year}{2020}) \bibinfo{pages}{89497--89509}.
\bibitem[{Huang et~al.(2019)Huang, Li, and Yao}]{huang2019survey}
\bibinfo{author}{C.~Huang}, \bibinfo{author}{Y.~Li}, \bibinfo{author}{X.~Yao},
\newblock \bibinfo{title}{A survey of automatic parameter tuning methods for metaheuristics},
\newblock \bibinfo{journal}{IEEE transactions on evolutionary computation} \bibinfo{volume}{24} (\bibinfo{year}{2019}) \bibinfo{pages}{201--216}.
\bibitem[{Lindauer et~al.(2022)Lindauer, Eggensperger, Feurer, Biedenkapp, Deng, Benjamins, Ruhkopf, Sass, and Hutter}]{smac3}
\bibinfo{author}{M.~Lindauer}, \bibinfo{author}{K.~Eggensperger}, \bibinfo{author}{M.~Feurer}, \bibinfo{author}{A.~Biedenkapp}, \bibinfo{author}{D.~Deng}, \bibinfo{author}{C.~Benjamins}, \bibinfo{author}{T.~Ruhkopf}, \bibinfo{author}{R.~Sass}, \bibinfo{author}{F.~Hutter},
\newblock \bibinfo{title}{Smac3: A versatile bayesian optimization package for hyperparameter optimization},
\newblock \bibinfo{journal}{Journal of Machine Learning Research} \bibinfo{volume}{23} (\bibinfo{year}{2022}) \bibinfo{pages}{1--9}.
\bibitem[{Shang et~al.(2020)Shang, Ishibuchi, He, and Pang}]{shang2020survey}
\bibinfo{author}{K.~Shang}, \bibinfo{author}{H.~Ishibuchi}, \bibinfo{author}{L.~He}, \bibinfo{author}{L.~M. Pang},
\newblock \bibinfo{title}{A survey on the hypervolume indicator in evolutionary multiobjective optimization},
\newblock \bibinfo{journal}{IEEE Transactions on Evolutionary Computation} \bibinfo{volume}{25} (\bibinfo{year}{2020}) \bibinfo{pages}{1--20}.
\bibitem[{Sahoo and Goswami(2023)}]{Sahoo_Goswami_2023}
\bibinfo{author}{S.~K. Sahoo}, \bibinfo{author}{S.~S. Goswami},
\newblock \bibinfo{title}{A comprehensive review of multiple criteria decision-making (mcdm) methods: Advancements, applications, and future directions},
\newblock \bibinfo{journal}{Decision Making Advances} \bibinfo{volume}{1} (\bibinfo{year}{2023}) \bibinfo{pages}{25–48}.
\bibitem[{Hwang and Yoon(2012)}]{hwang2012multiple}
\bibinfo{author}{C.~Hwang}, \bibinfo{author}{K.~Yoon}, \bibinfo{title}{Multiple Attribute Decision Making: Methods and Applications A State-of-the-Art Survey}, Lecture Notes in Economics and Mathematical Systems, \bibinfo{publisher}{Springer Berlin Heidelberg}, \bibinfo{year}{2012}.
\bibitem[{Vafaei et~al.(2018)Vafaei, Ribeiro, and Camarinha-Matos}]{vafaei2018data}
\bibinfo{author}{N.~Vafaei}, \bibinfo{author}{R.~A. Ribeiro}, \bibinfo{author}{L.~M. Camarinha-Matos},
\newblock \bibinfo{title}{Data normalisation techniques in decision making: case study with topsis method},
\newblock \bibinfo{journal}{International journal of information and decision sciences} \bibinfo{volume}{10} (\bibinfo{year}{2018}) \bibinfo{pages}{19--38}.
\bibitem[{{GEBCO Compilation Group}(2023)}]{gebco2023}
\bibinfo{author}{{GEBCO Compilation Group}}, \bibinfo{title}{Gebco 2023 grid (doi:10.5285/f98b053b-0cbc-6c23-e053-6c86abc0af7b)}, \bibinfo{year}{2023}.
\bibitem[{Jonsen et~al.(2005)Jonsen, Flemming, and Myers}]{jonsen2005robust}
\bibinfo{author}{I.~D. Jonsen}, \bibinfo{author}{J.~M. Flemming}, \bibinfo{author}{R.~A. Myers},
\newblock \bibinfo{title}{Robust state--space modeling of animal movement data},
\newblock \bibinfo{journal}{Ecology} \bibinfo{volume}{86} (\bibinfo{year}{2005}) \bibinfo{pages}{2874--2880}.
\bibitem[{Melo-Merino et~al.(2020)Melo-Merino, Reyes-Bonilla, and Lira-Noriega}]{melo2020ecological}
\bibinfo{author}{S.~M. Melo-Merino}, \bibinfo{author}{H.~Reyes-Bonilla}, \bibinfo{author}{A.~Lira-Noriega},
\newblock \bibinfo{title}{Ecological niche models and species distribution models in marine environments: A literature review and spatial analysis of evidence},
\newblock \bibinfo{journal}{Ecological Modelling} \bibinfo{volume}{415} (\bibinfo{year}{2020}) \bibinfo{pages}{108837}.
\bibitem[{Fan et~al.(2022)Fan, Yang, Yang, Wu, and Vladimir}]{fan2022review}
\bibinfo{author}{A.~Fan}, \bibinfo{author}{J.~Yang}, \bibinfo{author}{L.~Yang}, \bibinfo{author}{D.~Wu}, \bibinfo{author}{N.~Vladimir},
\newblock \bibinfo{title}{A review of ship fuel consumption models},
\newblock \bibinfo{journal}{Ocean engineering} \bibinfo{volume}{264} (\bibinfo{year}{2022}) \bibinfo{pages}{112405}.
\bibitem[{MacGillivray et~al.(2022)MacGillivray, Ainsworth, Zhao, Dolman, Hannay, Frouin-Mouy, Trounce, and White}]{MacGillivray2022}
\bibinfo{author}{A.~O. MacGillivray}, \bibinfo{author}{L.~M. Ainsworth}, \bibinfo{author}{J.~Zhao}, \bibinfo{author}{J.~N. Dolman}, \bibinfo{author}{D.~E. Hannay}, \bibinfo{author}{H.~Frouin-Mouy}, \bibinfo{author}{K.~B. Trounce}, \bibinfo{author}{D.~A. White},
\newblock \bibinfo{title}{{A functional regression analysis of vessel source level measurements from the Enhancing Cetacean Habitat and Observation (ECHO) database}},
\newblock \bibinfo{journal}{The Journal of the Acoustical Society of America} \bibinfo{volume}{152} (\bibinfo{year}{2022}) \bibinfo{pages}{1547--1563}.
\bibitem[{Deo et~al.(2024)Deo, Venkateshwaran, and Jaiman}]{deo2024continual}
\bibinfo{author}{I.~K. Deo}, \bibinfo{author}{A.~Venkateshwaran}, \bibinfo{author}{R.~K. Jaiman},
\newblock \bibinfo{title}{Continual learning of range-dependent transmission loss for underwater acoustic using conditional convolutional neural net},
\newblock \bibinfo{journal}{arXiv preprint arXiv:2404.08091}  (\bibinfo{year}{2024}).
\bibitem[{Deo and Jaiman(2022)}]{deo2022predicting}
\bibinfo{author}{I.~K. Deo}, \bibinfo{author}{R.~Jaiman},
\newblock \bibinfo{title}{Predicting waves in fluids with deep neural network},
\newblock \bibinfo{journal}{Physics of Fluids} \bibinfo{volume}{34} (\bibinfo{year}{2022}).
\bibitem[{Brissaud et~al.(2022)Brissaud, Näsholm, Turquet, and Le~Pichon}]{Brissaud2022}
\bibinfo{author}{Q.~Brissaud}, \bibinfo{author}{S.~P. Näsholm}, \bibinfo{author}{A.~Turquet}, \bibinfo{author}{A.~Le~Pichon},
\newblock \bibinfo{title}{{Predicting infrasound transmission loss using deep learning}},
\newblock \bibinfo{journal}{Geophysical Journal International} \bibinfo{volume}{232} (\bibinfo{year}{2022}) \bibinfo{pages}{274--286}.
\bibitem[{Audet and G(1974)}]{Audet1974AESDSP}
\bibinfo{author}{J.~J. Audet}, \bibinfo{author}{V.~G. G}, \bibinfo{title}{AESD Sound-Speed Profile Retrieval System (RSVP)}, \bibinfo{publisher}{Acoustic Environmental Support Detachment, Office of Naval Research}, \bibinfo{year}{1974}.

\end{thebibliography}





\end{document}